\documentclass{amsart}
\usepackage{amssymb,amsmath}

\theoremstyle{plain}
\newtheorem{thm}{Theorem}[section]
\newtheorem{theorem}[thm]{Theorem}
\newtheorem{lem}[thm]{Lemma}

\newtheorem{cor}[thm]{Corollary}
\newtheorem{prop}[thm]{Proposition}

\theoremstyle{definition}

\theoremstyle{remark}
\newtheorem{rem}[thm]{Remark}

\newtheorem*{remark*}{Remark}

\numberwithin{equation}{section}

\newcommand{\cF}{{\mathcal{F}}}

\newcommand{\cS}{{\mathcal{S}}}

        \newcommand{\field}[1]{{\mathbb{#1}}}
        \newcommand{\NN}{\field{N}}
        \newcommand{\ZZ}{\field{Z}}
        
        \newcommand{\RR}{\field{R}}

\allowdisplaybreaks

\begin{document}

\title[Schr\"odinger operators with magnetic wells]{Semiclassical
spectral asymptotics for a two-dimensional magnetic Schr\"odinger
operator II: The case of degenerate wells}

\author{Bernard Helffer}

\address{D\'epartement de Math\'ematiques, B\^atiment 425, Universit\'e
Paris-Sud 11, F-91405 Orsay C\'edex, France}
\email{Bernard.Helffer@math.u-psud.fr}

\author{Yuri A. Kordyukov }
\address{Institute of Mathematics, Russian Academy of Sciences, 112 Chernyshevsky
str. 450077 Ufa, Russia} \email{yurikor@matem.anrb.ru}

\thanks{Y.K. is partially supported by the Russian Foundation of Basic
Research (grant 09-01-00389).}

\begin{abstract}
We continue our study of a
magnetic Schr\"odinger operator on a two-dimensional compact
Riemannian manifold in the case when the minimal value of the module
of the magnetic field is strictly positive. We analyze the case when
the magnetic field has degenerate magnetic wells. The main result of
the paper is an asymptotics of the groundstate energy of the operator
in the semiclassical limit. The upper bounds are improved in the
case when we have a localization by a miniwell effect of lowest
order. These results are applied to prove the existence of an
arbitrary large number of spectral gaps in the semiclassical limit
in the corresponding periodic setting.
\end{abstract}

\subjclass[2000]{35P20, 35J10, 47F05, 81Q10}

\keywords{magnetic Schr\"odinger operator, eigenvalue asymptotics,
magnetic wells, semiclassical limit, spectral gaps}

\date{}
 \maketitle

\section{Preliminaries and main result}

Let $ M$ be a compact connected oriented manifold of dimension
$n\geq 2$ (possibly with boundary). Let $g$ be a Riemannian metric
and $\bf B$ a real-valued closed 2-form on $M$. Assume that $\bf B$
is exact and choose a real-valued 1-form $\bf A$ on $M$ such that
$d{\bf A} = \bf B$. Thus, one has a natural mapping
\[
u\mapsto ih\,du+{\bf A}u
\]
from $C^\infty_c(M)$ to the space $\Omega^1_c(M)$ of smooth,
compactly supported one-forms on $M$. The Riemannian metric allows
to define scalar products in these spaces and consider the adjoint
operator
\[
(ih\,d+{\bf A})^* : \Omega^1_c(M)\to C^\infty_c(M)\,.
\]
A Schr\"odinger operator with magnetic potential $\bf A$ is defined
by the formula
\[
H^h = (ih\,d+{\bf A})^* (ih\,d+{\bf A}).
\]
Here $h>0$ is a semiclassical parameter, which is assumed to be
small. If $M$ has non-empty boundary, we will assume that
 the metric $g$ and the magnetic potential $\bf A$ are
smooth up to boundary and the operator $H^h$ satisfies the
Dirichlet boundary conditions.

We are interested in semiclassical asymptotics of the low-lying
eigenvalues of the operator $H^h$. This problem was studied in
\cite{FoHel1,miniwells,HKI,HM,HM01,HelMo3,HelMo4,HelSj7,LuPa2,Mat,MatUeki,Mont,PiFeSt,Ray2,Ray3}
(see \cite{Fournais-Helffer:book,luminy} for surveys).

In this paper, we study the problem in a particular situation. We
suppose that $M$ is two-dimensional. Then we can write ${\bf B}=b\,
dx_g$, where $b\in C^\infty(M)$ and $dx_g$ is the Riemannian volume
form. Let
\[
b_0=\min_{x\in M} |b(x)|\,.
\]
We furthermore assume that:
\begin{itemize}
  \item $b_0>0$;
  \item the set $\{x\in M : |b(x)| =b_0\}$ is a smooth curve $\gamma$,
  which is contained in the interior of $M$;
\item there is a constant $C>0$ such that for all $x$
in some neighborhood of $\gamma$ the estimates hold:
\begin{equation}\label{YK:B1}
C^{-1}d(x,\gamma)^2\leq |b(x)|-b_0 \leq C
d(x,\gamma)^2\,.
\end{equation}
\end{itemize}

The main purpose is to give  an asymptotics of the groundstate energy
of the operator $H^h$. Denote by $N$ the external unit normal vector
to $\gamma$. Let $\tilde{N}$ denote the natural extension of $N$ to
a smooth normalized vector field on $M$, whose integral curves
starting from a point $x$ in  a tubular neighborhood of $\gamma$ are
the minimal geodesics  to $\gamma$. Consider the function $\beta_2$
on $\gamma$ given by
\begin{equation}\label{e:def-beta}
\beta_2(x)=\tilde{N}^2 |b(x)|\,, \quad x\in \gamma\,.
\end{equation}
By (\ref{YK:B1}), it is easy to see that
\[
\beta_2(x)>0\,, \quad x\in \gamma\,.
\]

\begin{theorem}~\\ \label{c:equiv-lambda0}
There exists $h_0>0$, such that, for any $h\in ]0,h_0]$,
\[
\lambda_0(H^h) =  hb_0 + h^2\,\frac{\mu_0}{4b_0} + \mathcal O ( h^{17/8}) \,.
\]
where
\[
\mu_0:=\inf_{x\in\gamma} \beta_2(x).
\]
\end{theorem}

The paper is organized as follows. In Section~\ref{s:upper} we
construct approximate eigenfunctions of the operator $H^h$. This
allows us to prove an upper bound in Theorem~\ref{c:equiv-lambda0}
(see Corollary \ref{Coroll}). In Section~\ref{s:lower} we prove a
lower bound for $\lambda_0(H^h)$ and complete the proof of
Theorem~\ref{c:equiv-lambda0}. In Section~\ref{s:miniwells},
assuming the existence of a unique non-degenerate miniwell for the
reduced spectral problem on $\gamma$, we construct more refined
approximate eigenfunctions of the operator $H^h$, improving the
upper bound of Section~\ref{s:upper}. In Section~\ref{s:gaps} we
consider the case when the magnetic field is periodic. We combine
the constructions of approximate eigenfunctions given in
Sections~\ref{s:upper} and~\ref{s:miniwells} with the results of
\cite{diff2006} to prove the existence of arbitrary large number of
gaps in the spectrum of the periodic operator $H^h$ in the
semiclassical limit.

{\bf Acknowledgements.} \\
We wish to thank the Erwin Schr\"odinger Institute in Vienna and the
organizers of the Conference ``QMath11 --- Mathematical Results in
Quantum Physics'' in Hradec Kralove for their hospitality and
support.

\section{Upper bounds}\label{s:upper}

\subsection{Approximate eigenfunctions: the main result} ~\\
The purpose of this section
is to construct approximate eigenfunctions for the operator $H^h$.
Denote by $R$ the scalar curvature of the Riemannian manifold
$(M,g)$.
\begin{theorem}\label{t:qmodes}~\\
For any $x\in \gamma $ and for any integer $k\geq 0$, there exist
$C$ and $h_0>0$, such that, for any $h\in ]0,h_0]$, there exists
$\Phi^h_k\in C^\infty_c(M), \Phi^h_k\neq 0\,,$ such that
\[
\left\|H^h\Phi^h_k-\lambda^h(k,x)\Phi^h_k\right\|\leq C h^{17/8}
\|\Phi^h_k\|\,,
\]
where
\begin{equation}\label{e:lambdakx}
\lambda^h(k,x)=(2k+1) hb_0 +
h^2\left[(2k^2+2k+1)\frac{\beta_2(x)}{4b_0}+\frac12
\left(k^2+k\right) R(x)\right]\,.
\end{equation}
\end{theorem}

When $k=0$, we get:
\begin{cor}\label{Coroll}~\\
For any $x\in \gamma $, there exist $C$ and $h_0>0$, such that, for
any $h\in ]0,h_0]$, there exists $\Phi^h_0\in C^\infty_c(M),
\Phi^h_0\neq 0\,,$ such that
\[
\left\|H^h\Phi^h_0-\lambda^h(x)\Phi^h_0\right\|\leq C h^{17/8}
\|\Phi^h_0\|\,,
\]
where
\[
\lambda^h(x)= hb_0 + h^2\,\frac{\beta_2(x)}{4b_0}\,.
\]
\end{cor}
This corollary  gives immediately the upperbound in Theorem
\ref{c:equiv-lambda0}.\\

The proof of Theorem \ref{t:qmodes}  is long, so we will split it in
different steps in the next subsections but let us previously discuss the interpretation of the coefficients.\\

\subsection{Geometrical interpretation of the coefficients}~\\
The term
\begin{equation}\label{e:Landau}
(2k+1) hb_0 + \frac12   h^2\left(k^2+k\right)R
\end{equation}
in the right-hand side of \eqref{e:lambdakx} has a natural
interpretation. It depends on whether $R$ is zero, positive or
negative  and in all three cases can be described in terms of
eigenvalues of the associated magnetic Laplacian with constant
magnetic field (Landau operator) on the corresponding simply
connected Riemann surface of constant curvature.

For $R=0$, it is a well-known fact that the eigenvalues of the
magnetic Laplacian with constant magnetic field $b$ on the flat
Euclidean plane
\[
H^h=\left(ih\frac{\partial}{\partial
x}-by\right)^2-h^2y^2\frac{\partial^2}{\partial y^2}
\]
are given by the Landau levels
\[
\lambda^h(k)=(2k+1) hb, \quad k\in \NN.
\]

In the case $R$ is negative, we consider the hyperbolic plane
$\mathbb H$, which is realized as the upper halfplane
\[
\mathbb H=\{(x,y)\in {\mathbb R}^2 : y>0\}\,,
\]
equipped with the Riemannian metric
\[
g=\frac{dx^2+dy^2}{y^2}\,.
\]
Then the scalar curvature is a negative constant $R=-2$. The
magnetic Laplacian with constant magnetic field $b$ on $\mathbb H$
is given by (see, for instance, \cite{Comtet87})
\[
H^h=y^2\left(ih\frac{\partial}{\partial
x}-by^{-1}\right)^2-h^2y^2\frac{\partial^2}{\partial y^2}\,.
\]
This operator first appeared in theory of automorphic forms, where
it is known as the Maas Laplacian. Its spectrum in $L^2(\mathbb H)$
was studied by Elstrodt in \cite{Elstrodt}. It consists of
absolutely continuous and discrete parts. The absolutely continuous
part is given by
\[
\sigma_{\rm{ac}}(H^h)=[h^2b^2+\frac14,+\infty[\,.
\]
The discrete part is empty if $0\leq hb\leq \frac12$ and, if $hb>
\frac12$, consists of a finite number of eigenvalues of infinite
multliplicity (hyperbolic Landau levels) given by
\[
\lambda^h(k)=(2k+1) hb - h^2\left(k^2+k\right)\,, \quad k\in\NN\,,
k<hb-\frac12\,.
\]
It is clear that the last formula agrees with \eqref{e:Landau}.

Finally, in the case $R$ is positive, we consider the
two-dimensional sphere
\[
S^2=\{(x,y,z)\in {\mathbb R}^3 : x^2+y^2+z^2=1\},
\]
equipped with the Riemannian metric induced by the standard
Euclidean metric in ${\mathbb R}^3$. Then the scalar curvature is a
positive constant $R=2$. In this context the magnetic Laplacian is
constructed as the Bochner Laplacian, acting on sections of a
Hermitian line bundle $\mathcal L$ with a compatible connection
$\nabla$ over $S^2$. The magnetic field is the curvature $2$-form
$\mathbf B$ of $\nabla$. It is constant if $\mathbf B$ is a scalar
multiple of the volume $2$-form $dx_g$:
\[
\mathbf B=s\,dx_g,\quad s\in {\mathbb R\,}\,.
\]

The construction of a Hermitian line bundle $\mathcal L$ with a
compatible connection $\nabla$ such that the associated magnetic
field is constant is a particular case of prequantization in
geometric quantization \cite{Kostant}. Such a line bundle exists if
and only if $s=n/2$ for some $n\in {\mathbb Z}$. For $s=n/2$, the
corresponding line bundle $({\mathcal L}_n,\nabla_n)$ can be
described as a complex line bundle associated with an
$S^1$-principal bundle $S^3\to S^2$ called the Hopf fibration and
the character $\chi_n : S^1\to S^1$ given by $ \chi_n(u)=u^n, u\in
S^1$. In physics literature, $({\mathcal L}_n,\nabla_n)$ is a
well-known Wu-Yang magnetic monopole \cite{Wu-Yang}, which provides
a natural topological interpretation of Dirac's monopole of magnetic
charge $g=nh/2e$.

The magnetic Laplacian $H_n$ acting on sections of ${\mathcal L}_n$
is defined as
\[
H_n=\nabla_n^*\nabla_n\,.
\]
Its spectrum was computed in \cite{Wu-Yang2} (see also
\cite{Kuwabara82,Kuwabara88}). It consists of a countable set of
eigenvalues (spherical Landau levels) given by
\begin{equation}\label{e:Landau-sph}
\frac{1}{2}|n|(2k+1)+k^2+k\,, \quad k\in {\mathbb N}\,,
\end{equation}
with multiplicity $|n|+2k+1$. The corresponding eigenvalues are
known as monopole harmonics. The formula \eqref{e:Landau-sph} agrees
with \eqref{e:Landau} if we take $h=1/n$ and
\[
H^h=\frac{1}{n^2}\nabla_n^*\nabla_n\,.
\]

Finally, let us mention the paper \cite{Ferapontov-Veselov}, which
states that the three types of magnetic Laplacians described above
are integrable in some sense and provides the complete description
of their spectra in the same way as it was done by Schr\"odinger for
the harmonic oscillator.

\subsection{Expanding operators in fractional powers of
$h$}\label{s:expand} ~\\
The approximate eigenfunctions $\Phi^h_k\in
C^\infty_c(M)$, which we are going to construct, will be supported
in a small neighborhood of the point $x\in \gamma$ appearing in the
theorem. We will consider some special local coordinate system with
coordinates $(s,t)$ in a neighborhood of $\gamma$ such that $\gamma$
corresponds to $t=0$. We will only apply our operator on functions
which are  products of cut-off functions with functions of the form
of linear combinations of terms like $h^\nu w(s,h^{-1/2}t)$ with $w$
in $C^\infty(S^1)\otimes \cS(\RR_t)$. These functions are
consequently $\mathcal O(h^\infty)$ outside a fixed neighborhood of
$\gamma$. We will start by doing the computations formally in the
sense that everything is determined modulo $\mathcal O(h^\infty)$,
and any smooth function will be replaced by its Taylor's expansion
at $\gamma$. It is then easy to construct non formal approximate
eigenfunctions.

Choose a normal coordinate system in a tubular neighborhood $U$ of
$\gamma$ with coordinates $X=(X_0,X_1)$ with $X_0=s$ and $X_1=t$.
Here $s\in [-L/2,L/2)\cong S^1_L=\RR/L\ZZ$ is the natural parameter
along $\gamma$ ($L$ is the length of $\gamma$), $\gamma$ is given by
the equation $t=0$, and $t\in (-\varepsilon_0,\varepsilon_0)$ is the
natural parameter along the geodesic, passing through the point on
$\gamma$ with the coordinates $(s,0)$ orthogonal to $\gamma$.

It is well known that in such coordinates the metric $g$ has the
form
\[
g=a(s,t)ds^2+dt^2,
\]
where
\[
a(s,0)=1\,.
\]
So we can write
\begin{equation}\label{e:expG}
g_{00}(s,t)=1+a_{1}(s)t+{a}_{2}(s)t^2+\mathcal O(t^3)\,, \quad
g_{01}(s,t) =0\,, \quad
 g_{11}(s,t) = 1\,.
\end{equation}
In particular, we have for the first coefficient of the inverse
matrix $(g^{ij})$
\begin{equation}\label{e:expG00}
g^{00}(s,t)=1-{a}_{1}(s)t-({a}_{2}(s)-{a}_{1}(s)^2)t^2+\mathcal
O(t^3).
\end{equation}
It is known (see \eqref{e:A1} and \eqref{e:A2} in the appendix) that
\[
a_{1}=-2\kappa, \quad a_2=-\frac12R+\kappa^2,
\]
where $\kappa$ is the mean curvature of $\gamma$.

Let us define $$
|g|=\det \{g_{ij}\}=g_{00}g_{11}-g_{01}^2\,,
$$
\[
{\bf A}=A_0ds+A_1dt\,,
\]
and
\begin{equation}\label{e:bfB}
{\bf B}=b(s,t)\sqrt{|g(s,t)|}ds\wedge dt, \quad
b(s,t)\sqrt{|g(s,t)|} =\frac{\partial A_1}{\partial
s}-\frac{\partial A_0}{\partial t}\,.
\end{equation}
The external unit normal vector to $\gamma$ has the form
\[
N=\frac{\partial }{\partial t},
\]
and its natural extension $\tilde{N}$ in the neighborhood of $\gamma$ is
\[
\tilde{N}=\frac{\partial }{\partial t}\,.
\]
Without loss of generality, we can assume that $b(x)>0$ for any
$x\in M$. Thus, we have
\begin{equation}\label{e:b}
b(s,0)=b_0,\quad \frac{\partial b}{\partial t}(s,0)=0, \quad
\beta_2(s)=\frac{\partial^2 b}{\partial t^2}(s,0)>0.
\end{equation}

We can assume that (after a gauge transformation)
\[
A_1(s,t)\equiv 0.
\]
Using \eqref{e:expG}, \eqref{e:bfB} and \eqref{e:b}, we obtain that
\begin{multline}\label{e:A0}
A_0(s,t) \\=\alpha_1(s)-b_0t-\frac14{a}_{1}(s)b_0t^2
-\frac16\left(\beta_2(s)+b_0({a}_{2}(s)-\frac14
{a}_{1}(s)^2)\right)t^3\\ +\mathcal O(t^4)\,,
\end{multline}
as $t \to 0\,$.
Our constructions will be local, in a neighborhood of $x\in \gamma$.
So, using gauge invariance, we can take $\alpha_1=0$.

We have
\[
H^h=\frac{1}{\sqrt{|g(X)|}}\sum_{0\leq \alpha,\beta\leq
1}\nabla^h_\alpha \left(\sqrt{|g(X)|} g^{\alpha\beta}(X)
\nabla^h_\beta\right)\,,
\]
where
\[
\nabla^h_\alpha= i h \frac{\partial}{\partial X_\alpha}+A_\alpha(X)\,,
\quad \alpha =0,1\,,
\]
or, taking into account \eqref{e:expG},
\begin{equation}\label{e:Hh}
H^h= g^{00}(s,t) (\nabla^h_0)^2 - h^2 \frac{{\partial}^2}{\partial
t^2} + i h\Gamma^0 \nabla^h_0  - h^2\Gamma^1
\frac{\partial}{\partial t}\,.
\end{equation}
where
\[
\Gamma^\alpha= \frac{1}{\sqrt{|g(X)|}}\sum_{0\leq \beta\leq
1}\frac{\partial}{\partial X_\beta} \left(\sqrt{|g(X)|}
g^{\beta\alpha}(X)\right)\,, \quad \alpha =0,1\,.
\]
By \eqref{e:expG} and \eqref{e:expG00}, it follows that
\begin{align*}
\Gamma^0=&\frac{1}{\sqrt{|g(X)|}}\frac{\partial}{\partial s}
\left(\sqrt{|g(X)|} g^{00}(X)\right) =
-\frac12a^\prime_{1}(s)t+\mathcal O(t^2)\,,
\\ \Gamma^1 = & \frac{1}{\sqrt{|g(X)|}}\frac{\partial}{\partial t}
\left(\sqrt{|g(X)|}\right)=
\frac12{a}_{1}(s)+\left({a}_{2}(s)-\frac12
{a}_{1}(s)^2\right)t+\mathcal O(t^2)\,.
\end{align*}

We now move the operator $H^h$ into the Hilbert space
$L^2(U,ds\,dt)$, considering the operator
\begin{equation}\label{e:def-hatH}
\hat{H}^h : =|g(X)|^{1/4}H_h |g(X)|^{-1/4}\,.
\end{equation}
By \eqref{e:Hh}, we get
\begin{equation}\label{e:form-hatH}
\hat{H}^h = g^{00}(s,t) (\hat\nabla^h_0)^2 + (\hat\nabla^h_1)^2 + i
h\Gamma^0 \hat\nabla^h_0 +i h\Gamma^1 \hat\nabla^h_1\,,
\end{equation}
where
\begin{equation}\label{e:defNabla}
\begin{aligned}
\hat{\nabla}^h_0& : =|g(X)|^{1/4} {\nabla}^h_0 |g(X)|^{-1/4}\\ &  =
{\nabla}^h_0-\frac14ih a^\prime_1(s)t+\mathcal O(t^2)\,,
\\ \hat{\nabla}^h_1& : =ih |g(X)|^{1/4} \frac{\partial}{\partial t}
|g(X)|^{-1/4}\\ &  = ih \left(\frac{\partial}{\partial
t}-\frac14a_1(s)-\frac12(a_2(s)-\frac12a_1(s)^2)t+\mathcal
O(t^2)\right)\,.
\end{aligned}
\end{equation}

Now we use the scaling $t=h^{1/2}t_1$ and expand the operator
$\hat{H}^h$ in power of $h^{1/2}$. By straightforward computations,
we obtain that
\[
\hat{H}^h=hH^{h,0}\,,
\]
where
\[
H^{h,0}=P_0+h^{1/2}P_1+hP_2+h^{3/2}{\mathcal R}_3(h)\,,
\]
with
\begin{align*}
P_0=& - \frac{{\partial}^2}{\partial t_1^2} + b^2_0t^2_1,\\
P_1=& - 2b_0t_1\left( i \frac{\partial}{\partial
s}+\frac14a_1(s)b_0t^2_1\right), \\
P_2=& -\frac{\partial^2}{\partial s^2}+\frac{3i}{2}a_1(s)
\frac{\partial}{\partial s} b_0t^2_1+\frac13b_0t^4_1\beta_2(s) \\
& +\frac{3i}{4}a^\prime_1(s)b_0t^2_1
+\frac{23}{48}a_1(s)^2b^2_0t^4_1-\frac{2}{3}a_2(s)b_0^2t^4_1
\\ & +\frac12 \left(a_2(s)-\frac38
a_1(s)^2\right),
\end{align*}
and ${\mathcal R}_3(h)$ a second order differential operator whose
coefficients are formal power series of the form $\sum_{j=0}^\infty
h^{j/2}a_{j}(s,t_1)$.

\subsection{Construction of approximate eigenfunctions}~\\
Let us fix
$x \in \gamma$ and $k\in {\mathbb N}$. Without loss of
generality, we can assume that $x$ corresponds to $0\in S^1_L$. We
introduce a self-adjoint second order differential operator
$Q^{h,0}$ in $L^2(S^1_L\times \RR)$, which is close to $H^{h,0}$
near $\{0\}\times \RR$, and construct approximate eigenfunctions
$\varphi_h$ of the operator $Q^{h,0}$ in the form
\[
\varphi_h=\varphi_0+h^{1/2}\varphi_1+h\varphi_2\,,
\]
such that
\[
\|({Q}^{h,0}-\lambda(h))\varphi_h\|=\mathcal O(h^{3/2})\,,
\]
with $\lambda(h)$ of the form
$\lambda(h)=\lambda_0+h^{1/2}\lambda_1+h\lambda_2$\,. (Observe that
such approximate eigenfunctions don't exist for the operator
$H^{h,0}$.) The fact that the operator $H^{h,0}$ is close to
$Q^{h,0}$ near $0$ will allow us to construct approximate
eigenfunctions for $H^{h,0}$, using appropriate $h$-dependent
cut-off functions, and complete the proof of Theorem~\ref{t:qmodes}.

The operator $Q^{h,0}$ is defined by the formula
\[
Q^{h,0}:=L_0+h^{1/2}L_1+hL_2\,,
\]
where
\begin{align*}
L_0=& - \frac{{\partial}^2}{\partial t_1^2} + b^2_0t^2_1\,,\\
L_1=& - 2b_0t_1\left( i \frac{\partial}{\partial
s}+\frac14a_1(s)b_0t^2_1\right)\,, \\
L_2=& -\frac{\partial^2}{\partial s^2}+\frac{3i}{2}a_1(s)
\frac{\partial}{\partial s} b_0t^2_1+\frac13b_0t^4_1\beta_2(0) \\
& -\frac{2}{3}a_2(0)b_0^2t^4_1 +\frac{1}{6}a_1(0)^2b^2_0t^4_1
+\frac12 \left(a_2(0)-\frac14 a_1(0)^2\right) \\ & +
\frac{3i}{4}a^\prime_1(s)b_0t^2_1 +\frac{5}{16}a_1(s)^2b^2_0t^4_1
-\frac{1}{16}a_1(s)^2\,.
\end{align*}

Now we write formal expansions in powers of $h^{1/2}$:
\[
\varphi_h\sim \sum_{j=0}^\infty h^{j/2}\varphi_j\,, \quad
\lambda(h)\sim \sum_{j=0}^\infty h^{j/2}\lambda_j\,,
\]
and express the cancellation of the coefficients of $h^{j/2}$ in the
formal expansion for $ (Q^{h,0}-\lambda(h))\varphi_h$
for $j=0,1,2\,$.\\
{\bf Step 1}~\\
At the first step we get
\begin{equation}\label{e:lambda_0}
L_0\varphi_0= \lambda_0\varphi_0\,.
\end{equation}
We have
\[
{\rm Sp}(L_0)=\{\mu_{m}=(2m+1)b_0: m\in \NN\}\,.
\]
The eigenfunction of $L_0$ associated to the eigenvalue $\mu_{m}$ is
\[
\psi_{m}(t_1)=\pi^{-1/4}b_0^{1/2} H_{m}(b_0^{1/2}t_1)
e^{-b_0t^2_1/2}\,,
\]
where $H_m$ is the Hermite polynomial:
\[
H_m(x)=(-1)^me^{x^2}\frac{d^m}{dx^m}(e^{-x^2})\,.
\]
The norm of $\psi_{m}$ in $L^2(\RR,dx)$ equals the norm of $H_m$ in
$L^2(\RR,e^{-x^2}\,dx)$\,, which is given by
\[
\|H_m\|=\sqrt{2^mm!\sqrt{m}}\,.
\]
For the given  fixed integer $k\geq 0$, we take a solution of
\eqref{e:lambda_0} in the form
\[
\lambda_0=(2k+1)b_0\,,\quad \varphi_0(s,t_1)=\chi_0(s)\psi_k(t_1)\,.
\]
We recall that the Hermite polynomials satisfy
\[
H_{m+1}(x)=2xH_m(x)-2mH_{m-1}(x)\,,
\]
and
\[
H^\prime_m(x)=2xH_m(x)-H_{m+1}\,.
\]
We also have
\begin{align*}
2xH_m= & H_{m+1}+2mH_{m-1},\\
4x^2H_m= & H_{m+2}+(4m+2)H_m+4m(m-1)H_{m-2}\,,\\
8x^3H_m= & H_{m+3}+(6m+6)H_{m+1}+12m^2H_{m-1}+8m(m-1)(m-2)H_{m-3}\,,\\
16x^4H_m(x)= & H_{m+4}(x)+(8m+12)H_{m+2}(x)+12(2m^2+2m+1)H_m(x)\\
&+16(2m^2-3m+1)mH_{m-2}(x)\\ & +16m(m-1)(m-2)(m-3)H_{m-4}(x)\,.
\end{align*}
Using these identities and the orthogonality of the Hermite polynomials,
we get
\begin{gather*}
\frac{1}{\|H_k\|^2}\langle t_1\psi_{k-1}\, , \psi_{k}\rangle =
\frac{1}{2b_0^{1/2}}\,,\quad
\frac{1}{\|H_k\|^2}\langle t_1\psi_{k+1}\,, \psi_{k}\rangle =  \frac{k+1}{b_0^{1/2}}\,,\\
\frac{1}{\|H_k\|^2}\langle t^3_1\psi_{k-3}, \psi_{k}\rangle =
\frac{1}{8b_0^{3/2}}\,,\quad
\frac{1}{\|H_k\|^2}\langle t^3_1\psi_{k-1}, \psi_{k}\rangle =  \frac{3}{4b_0^{3/2}}k\,,\\
\frac{1}{\|H_k\|^2}\langle t^3_1\psi_{k+3}, \psi_{k}\rangle =
\frac{1}{b_0^{3/2}}(k+3)(k+2)(k+1)\,,
\frac{1}{\|H_k\|^2}\langle t^3_1\psi_{k+1}, \psi_{k}\rangle =  \frac{3}{2b_0^{3/2}}(k+1)^2\,,\\
\frac{1}{\|H_k\|^2}\langle t^2_1\psi_{k}, \psi_{k}\rangle =
\frac{1}{2b_0}(2k+1)\,,\quad \frac{1}{\|H_k\|^2}\langle t^4_1
\psi_{k}, \psi_{k}\rangle = \frac{3}{4b_0^2}(2k^2+2k+1)\,.
\end{gather*}

We will also use several times the following identities:
\[
\frac{1}{\|H_k\|^2} \langle
\left(t^3_2D_{t_2}-\frac{3}{2}it^2_2\right)\psi_k, \psi_k\rangle=0\,,
\]
\begin{multline*}
\frac{1}{\|H_k\|^2}\langle \frac{1}{6}t_2^3\psi_{k+3} +
(3k+3)t_2^3\psi_{k+1}
-6k^2t_2^3\psi_{k-1}-\frac{4}{3}k(k-1)(k-2)t_2^3\psi_{k-3}, \psi_k\rangle\\
= \frac{1}{b^{3/2}_0}\left(15k^2+15k+ \frac{11}{2}\right)\,,
\end{multline*}
and
\begin{multline*}
\frac{1}{\|H_k\|^2} \langle
t_2\Big(\frac{1}{6}\psi_{k+3}+(3k+3)\psi_{k+1}
-6k^2\psi_{k-1}-\frac{4}{3}k(k-1)(k-2)\psi_{k-3}\Big),\psi_k\rangle\\
= \frac{3}{b_0^{1/2}}(2k+1)\,.
\end{multline*}
{\bf Step 2}\\
At the second step, we obtain
\begin{equation} \label{e:lambda_1}
(L_0-\lambda_0)\varphi_1=\lambda_1\varphi_0-L_1\varphi_0\,,
\end{equation}
or, more explicitly,
\begin{multline*}
\left(-\frac{\partial^2 }{\partial t_1^2}
+b^2_0t^2_1-(2k+1)b_0\right)\varphi_1\\
\begin{aligned}
= & \lambda_1\chi_0(s)\psi_k(t_1)
+ i\chi^\prime_0(s) b^{1/2}_0(\psi_{k+1}(t_1)+2k\psi_{k-1}(t_1))\\
& + a_1(s)\chi_0(s) b^{1/2}_0
\Big(\frac{1}{16}\psi_{k+3}(t_1)+\left(\frac{3}{8}k+\frac38\right)\psi_{k+1}(t_1)\\
& +\frac{3}{4}k^2\psi_{k-1}(t_1)
+\frac{1}{2}k(k-1)(k-2)\psi_{k-3}(t_1)\Big)\,.
\end{aligned}
\end{multline*}
The equation \eqref{e:lambda_1} has a solution only when
\[
\lambda_1=0\,,
\]
and this solution can be taken in the form
\begin{align*}
\varphi_1(s,t_1)= & \frac{i}{2b^{1/2}_0}\chi^\prime_0(s)[\psi_{k+1}(t_1)-2k\psi_{k-1}(t_1)]\\
& + a_1(s)\chi_0(s)\frac{1}{2b^{1/2}_0}
\Big(\frac{1}{48}\psi_{k+3}(t_1)+\left(\frac{3}{8}k+\frac38\right)\psi_{k+1}(t_1)\\
& - \frac{3}{4}k^2\psi_{k-1}(t_1)
-\frac{1}{6}k(k-1)(k-2)\psi_{k-3}(t_1)\Big)\,.
\end{align*}
{\bf Step 3}\\
At the third step, we have
\begin{equation} \label{e:lambda_2}
(L_0-\lambda_0)\varphi_2=\lambda_2\varphi_0-L_2\varphi_0-L_1\varphi_1\,.
\end{equation}
This equation has a solution if and only if, for any $s$, its
right-hand side is orthogonal to $\psi_k(t_1)$.

First, let us compute $L_2\varphi_0\,$:
\begin{align*}
L_2\varphi_0=& -\chi^{\prime\prime}_0(s)\psi_k(t_1)
+\frac{3i}{2}a_1(s)\chi^\prime_0(s) b_0t^2_1\psi_k(t_1)+\frac13\beta_2(0)\chi_0(s)b_0t^4_1\psi_k(t_1) \\
& -\frac{2}{3}a_2(0)\chi_0(s)b_0^2t^4_1\psi_k(t_1)
+\frac{1}{6}a_1(0)^2\chi_0(s) b^2_0t^4_1\psi_k(t_1)\\ & +\frac12
\left(a_2(0)-\frac14 a_1(0)^2\right)\chi_0(s)\psi_k(t_1) +
\frac{3i}{4}a^\prime_1(s)\chi_0(s) b_0t^2_1\psi_k(t_1)\\
& +\frac{5}{16}a_1(s)^2\chi_0(s) b^2_0t^4_1\psi_k(t_1)
-\frac{1}{16}a_1(s)^2\chi_0(s)\psi_k(t_1)\,.
\end{align*}
Multiplying by $\psi_k$ and integrating with respect to $t_1$ gives:
\begin{multline*}
\frac{1}{\|H_k\|^2}\int L_2\varphi_0 (s,t_1)\psi_k(t_1)\,dt_1\\
\begin{aligned}
=& -\chi^{\prime\prime}_0(s)
+ a_1(s)\chi^\prime_0(s)i\left(\frac{3}{2}k+\frac{3}{4}\right)+\frac{1}{4b_0}(2k^2+2k+1)\beta_2(0)\chi_0(s) \\
& -(k^2+k) a_2(0)\chi_0(s) +\frac{1}{4}(k^2+k)a_1(0)^2\chi_0(s)\\
& + \frac{3i}{8}(2k+1)a^\prime_1(s)\chi_0(s)
+\left(\frac{15}{32}k^2+\frac{15}{32}k+\frac{11}{64}\right)
a_1(s)^2\chi_0(s)\,.
\end{aligned}
\end{multline*}
Next, let us compute $L_1\varphi_1$:
\begin{align*}
L_1\varphi_1=&
b^{1/2}_0\chi^{\prime\prime}_0(s)[t_1\psi_{k+1}(t_1)-2kt_1\psi_{k-1}(t_1)]\\
& - a^{\prime}_1(s)\chi_0(s)b^{1/2}_0i
(\frac{1}{48}t_1\psi_{k+3}(t_1)+\left(\frac{3}{8}k+\frac38\right)t_1\psi_{k+1}(t_1)\\
& -\frac{3}{4}k^2t_1\psi_{k-1}(t_1)
-\frac{1}{6}k(k-1)(k-2)t_1\psi_{k-3}(t_1))\\
& - a_1(s)\chi^{\prime}_0(s)b^{1/2}_0i
\Big(\frac{1}{48}t_1\psi_{k+3}(t_1)+\left(\frac{3}{8}k+\frac38\right)t_1\psi_{k+1}(t_1)\\
& -\frac{3}{4}k^2t_1\psi_{k-1}(t_1)
-\frac{1}{6}k(k-1)(k-2)t_1\psi_{k-3}(t_1)\Big)\\
&
-\frac{i}{4}a_1(s)\chi^\prime_0(s)b^{3/2}_0[t^3_1\psi_{k+1}(t_1)-2kt^3_1\psi_{k-1}(t_1)]\\
& - \frac14a_1(s)^2\chi_0(s)b^{3/2}_0
\Big(\frac{1}{48}t^3_1\psi_{k+3}(t_1)+\left(\frac{3}{8}k+\frac38\right)t^3_1\psi_{k+1}(t_1)\\
& -\frac{3}{4}k^2t^3_1\psi_{k-1}(t_1)
-\frac{1}{6}k(k-1)(k-2)t^3_1\psi_{k-3}(t_1)\Big)\,.
\end{align*}
Multiplying by $\psi_k$ and integrating with respect to $t_1$ gives:
\begin{multline*}
\frac{1}{\|H_k\|^2}\int L_1\varphi_1 (s,t_1)\psi_k(t_1)\,dt_1\\
\begin{aligned}
=& \chi^{\prime\prime}_0(s)-
a^{\prime}_1(s)\chi_0(s)i\left(\frac{3}{4}k+\frac38\right)-
a_1(s)\chi^{\prime}_0(s)i\left(\frac{3}{2}k+\frac34\right)\\
& - \frac14a_1(s)^2\chi_0(s) \left(\frac{15}{8}k^2+\frac{15}{8}k
+\frac{11}{16}\right)\,.
\end{aligned}
\end{multline*}
Thus, multiplying the right-hand side of \eqref{e:lambda_2} by
$\psi_k$ and integrating with respect to $t_1$ we obtain
\begin{multline*}
\frac{1}{\|H_k\|^2}\langle \lambda_2\varphi_0-L_2\varphi_0-L_1\varphi_1, \psi_k\rangle\\
\begin{aligned}
=& \lambda_2\chi_0(s)-\frac{1}{\|H_k\|^2}\int L_2\varphi_0
(s,t_1)\psi_k(t_1)\,dt_1 -\frac{1}{\|H_k\|^2}\int L_1\varphi_1 (s,t_1)\psi_k(t_1)\,dt_1\\
= & \lambda_2\chi_0(s)- \frac{\beta_2(0)}{4b_0}(2k^2+2k+1)\chi_0(s)+
\left(k^2+k\right)\left(a_2(0)-\frac14a_1(0)^2\right)\chi_0(s)\,.
\end{aligned}
\end{multline*}
So the solvability condition for \eqref{e:lambda_2} holds if we put
\[
\lambda_2=\frac{\beta_2(0)}{4b_0}(2k^2+2k+1)-
\left(k^2+k\right)\left(a_2(0)-\frac14a_1(0)^2\right)\,,
\]
and, in this case, there exists a solution $\varphi_2$ of
\eqref{e:lambda_2}.

Thus, we arrive at the following proposition.

\begin{prop}~\\
For any $\chi_0\in C^\infty_c(\mathbb R)$, there exists an
approximate eigenfunction $\varphi_h$ of the operator $Q^{h,0}$ in
the form
\begin{equation}\label{e:varphih}
\varphi_h=\varphi_0+h^{1/2}\varphi_1+h\varphi_2\,,
\end{equation}
where $\varphi_0$ is given by
\begin{equation}\label{e:varphi_0}
\varphi_0(s,t_1)=\chi_0(s)\psi_k(t_1)\,,
\end{equation}
$\varphi_1$ and $\varphi_2$ have the form
\begin{align}\label{e:varphi_1}
\varphi_1(s,t_1)&=A_0(s,t_1)\chi_0(s)+A_1(s,t_1)\chi_0^\prime(s)\,,\\
\label{e:varphi_2}
\varphi_2(s,t_1)&=B_0(s,t_1)\chi_0(s)+B_1(s,t_1)\chi_0^\prime(s)+B_2(s,t_1)\chi_0^{\prime\prime}(s)\,,
\end{align}
where $A_0$, $A_1$, $B_0$, $B_1$ and $B_2$ are some smooth
functions, such that
\[
\|({Q}^{h,0}-\lambda(h))\varphi_h\|=\mathcal O(h^{3/2})\,,
\]
with
\[
\lambda(h)=\lambda_0+h\lambda_2\,,
\]
where
\[
\lambda_0=(2k+1)b_0\,,\quad
\lambda_2=\frac{\beta_2(0)}{4b_0}(2k^2+2k+1)\,-
\left(k^2+k\right)\left(a_2(0)-\frac14a_1(0)^2\right)\,.
\]
\end{prop}
\noindent Observe that $\lambda(h)=h^{-1}\lambda^h(k,x)$, where
$\lambda^h(k,x)$ is given by \eqref{e:lambdakx}.\\

Now we consider the  function $\Phi_h$, as  constructed in \eqref{e:varphih}  in the previous proposition but with
$\chi_0$, depending on $h$:
\[
\chi_0(s):=E_h(s)=c h^{-\beta/2}\chi(s) \exp
\left(-\frac{s^2}{2h^{2\beta}} \right)\,.
\]
Here $\beta\in ]0,1/2[$ is some constant, which we will choose
later, and $\chi\in C^\infty_c(\mathbb R)$ is a cut-off function,
which equals $1$ in a neighborhood of $0$.

Observe that
\begin{equation}\label{YK:2}
\|s^m E_h(s)\| = \left(h^{-\beta} \int_{\RR} s^{2m}\exp
\left(-\frac{s^2}{h^{2\beta}} \right) ds \right)^{1/2}\\ =C_1
h^{\beta m}\,,
\end{equation}
and, furthermore,
\begin{equation}\label{YK:3}
\|s^m E_h^\prime(s)\| = C_2h^{\beta (m-1)}, \quad \|s^m
E^{\prime\prime}_h(s)\| = C_3h^{\beta (m-2)}\,.
\end{equation}
We have
\[
({Q}^{h,0}-\lambda(h))\Phi_h = h^{3/2}(L_2\varphi_1+L_1\varphi_2
-\lambda_2\varphi_1)+h^{2}(L_2\varphi_2 -\lambda_2\varphi_2)\,.
\]
Since $L_1$ is a first order differential operator in $s$ and $L_2$
is a second order differential operator in $s$, by
\eqref{e:varphi_1} and \eqref{e:varphi_2}, it follows that
\[
L_2\varphi_1+L_1\varphi_2=F_0E_h+F_1E_h^\prime+F_2
E_h^{\prime\prime}+F_3E_h^{\prime\prime\prime}\,,
\]
where $F_j, j=0,1,2,3$, are some functions. Using \eqref{YK:2} and
\eqref{YK:3}, we obtain the existence of constants $C>0$ and $h_0>0$ such  that
\[
\|({Q}^{h,0}-\lambda(h))\Phi_h\|\leq C h^{3/2-3\beta} \|\Phi_h\|\,,\, \forall h\in]0,h_0]\,.
\]
Then  we have
\begin{multline*}
H^{h,0}-Q^{h,0} = \\
\begin{aligned}
 & h\Bigg\{\frac13 b_0 t_1^4(\beta_2(s)-\beta_2(0))
-\frac{2}{3}(a_2(s)-a_2(0))b_0^2t^4_1
+\frac{1}{6}(a_1(s)^2-a_1(0)^2)b^2_0t^4_1\\ & +\frac12\left[
\left(a_2(s)-\frac14 a_1(s)^2\right)-\left(a_2(0)-\frac14
a_1(0)^2\right)\right]\Bigg\}+h^{3/2}{\mathcal R}_3(h)\,,
\end{aligned}
\end{multline*}
which immediately implies that
\[
\|(H^{h,0}-Q^{h,0})\Phi_h\|\leq C h^{\gamma_0(\beta)}\|\Phi_h\|\,,
\]
with
\[
\gamma_0(\beta)=\min(1+\beta,3/2-2\beta)\,,
\]
and
\[
\|(H^{h,0}-\lambda(h))\Phi_h\|\leq C h^{\gamma(\beta)}\|\Phi_h\|\,,
\]
with
\[
\gamma(\beta)=\min(1+\beta,3/2-3\beta)\,.
\]
Setting
\[
\beta=\frac18\,,
\]
we obtain the existence of a constant $C$ such  that
\[
\|(H^{h,0}-\lambda(h))\Phi_h\|\leq C h^{9/8} \|\Phi_h\|\,.
\]

\section{Lower bounds}\label{s:lower}
In this section, we will prove the lower bound for the groundstate
energy $\lambda_0(H^h)$ of the operator $H^h$. As above, we write
\[
\mu_0:=\inf_{s\in \gamma}\beta_2(s)\,.
\]
\begin{thm}~\\
There exist $C$ and $h_0>0$, such that for any $h\in ]0,h_0]$
\[
\lambda_0(H^h)\geq hb_0 + h^2\,\frac{\mu_0}{4b_0}-C h^{19/8}\,.
\]
\end{thm}

First, we recall a general lower bound due to Montgomery
\cite{Mont}. Suppose that $U$ is a domain in $M$. Then, for any
$u\in C^\infty_c(U)$, the following estimate holds:
\begin{equation}\label{e:Mont}
\|(ih\,d+{\bf A})u\|^2_U\geq h \left|\int_U b |u|^2dx_g\right|\,.
\end{equation}
This fact is an immediate consequence of a Weitzenb\"ock-Bochner
type identity.

From \eqref{e:Mont}, it follows that we can restrict our
considerations to  any sufficiently small neighborhood $\Omega$ of
$\gamma$. We denote by $H^h_{D}$ the Dirichlet realization of the
operator $H^h$ in $L^2(\Omega,dx_g)\,$.

The estimate \eqref{e:Mont} implies that
\[
\tau h H^h_{D} +(1-\tau h)hb\leq H^h_{D}, \quad 0<\tau<h^{-1}\,.
\]
Taking $\tau=h^{-1/2}\,$, we obtain
\[
h^{1/2} (H^h_{D}-hb +h^{1/2}b)\leq H^h_{D}\,, \quad 0<h<1\,.
\]
Consider the Dirichlet realization $P^h_D$ of the operator
\begin{equation}
P^h=H^h-hb+h^{1/2}(b-b_0)\,,
\end{equation}
in $L^2(\Omega,dx_g)$. Then we have
\begin{equation}\label{e:lambdaj}
hb_0+h^{1/2} \lambda_0(P^h_{D})\leq \lambda_0(H^h_{D})\,.
\end{equation}
Therefore, the desired lower bound for $\lambda_0(H^h_{D})$ is an
immediate consequence of the following proposition.
\begin{prop}\label{t:PAV}~\\
There exist $C_0$ and $h_0>0$ such that
\[
\lambda_0(P^h_D)\geq h^{3/2}\,\frac{\mu_0}{4b_0} - C_0h^{15/8}\,,
\quad h\in ]0,h_0]\,.
\]
\end{prop}

To prove Proposition~\ref{t:PAV}, we will follow the lines of the proof
of \cite[Theorem 7.4]{HM01}. First, we observe that the upper bound
in Theorem~\ref{c:equiv-lambda0} and \eqref{e:lambdaj} imply an
upper bound for $\lambda_0(P^h_D)$:
\begin{equation}\label{e:upperPh}
\lambda_0(P^h_D)\leq h^{3/2}\,\frac{\mu_0}{4b_0} + Ch^{13/8}\,.
\end{equation}

Denote by $u_h$ an eigenfunction of $P^h_D$ associated with
$\lambda_0(P^h_D)$:
\begin{equation}\label{e:Vhj}
P^h_Du_h= \lambda_0(P^h_D)u_h\,.
\end{equation}
By a straightforward repetition of the arguments of \cite{HM01}, we
can easily show the following analogue of Lemmas 7.10 and 7.11 in
\cite{HM01}\footnote{There are a few inaccuracies in \cite{HM01},
concerning Lemma 7.11. For the erratum, see
http://www.math.u-psud.fr/$\sim$helffer/erratum164II.pdf }.

\begin{lem}\label{l:weight}~\\
For any real $k\geq 0$, we have
\[
\| d(x,\gamma)^ku_h\|_{L^2(\Omega,dx_g)}\leq
C_{k,j}(h^{k/2}+h^{(3k+1)/8}) \| u_h\|_{L^2(\Omega,dx_g)}\,.
\]
For any $\alpha =0,1$ and any $k\geq 0$, we have
\[
\| d(x,\gamma)^k\nabla^h_\alpha u_h\|_{L^2(\Omega,dx_g)}\leq
C_{k,j}(h^{(k+1)/2}+h^{(3k+5)/8})\|u_h\|_{L^2(\Omega,dx_g)}\,.
\]
\end{lem}

As above, take normal local coordinates $X=(X_0,X_1)=(s,t)$ near
$\gamma$ such that $\gamma$ corresponds to $t=0$ and assume that
they are defined on $\Omega$. Thus, we have
\[
b(X)=b_0+\frac12\beta_2(s) t^2+\mathcal O(t^3)\,.
\]
and, for the metric coefficients,
\[
g_{00}(s,t)=1+a_{1}(s)t+{a}_{2}(s)t^2+\mathcal O(t^3)\,, \quad
g_{01}(s,t) =0\,, \quad
 g_{11}(s,t) = 1\,.
 \]
We can choose a magnetic potential $A$ such that
\[
A_0(s,t)=\alpha_1(s)-b_0t-\frac14{a}_{1}(s)b_0t^2+\mathcal O(t^4),
\quad A_1(X)=0\,.
\]
We have
\[
P^h_{D}=\sum_{0\leq \alpha,\beta\leq 1} \frac{1}{\sqrt{|g(X)|}}
\nabla^h_\alpha \sqrt{|g(X)|} g^{\alpha\beta}(X)
\nabla^h_\beta-hb(X)+h^{1/2} (b(X)-b_0)\,,
\]
so its quadratic form is given by
\begin{multline*}
(P^h_{D}u,u)= \int_\Omega \sum_{0\leq \alpha,\beta\leq 1}
g^{\alpha\beta}(X) \nabla^h_\alpha u(X) \overline{\nabla^h_\beta
u(X)} \sqrt{g(X)}dX \\ -h\int_\Omega b(X)|u(X)|^2\sqrt{g(X)}\,dX
+h^{1/2} \int_\Omega (b(X)-b_0)|u(X)|^2\sqrt{g(X)}\,dX\,.
\end{multline*}
Note that
\[
P^h_{D}\geq 0\,.
\]

Now we move the operator $P_D^{h}$ into the Hilbert space
$L^2(\Omega,dX)$, using the unitary change of variables
$v=|g(X)|^{1/4}u$. For the corresponding operator
$$ \hat{P}_D^{h}=|g(X)|^{1/4}P_D^{h} |g(X)|^{-1/4}$$ in
$L^2(\Omega,dX)\,$, we obtain
\begin{multline*}
(\hat{P}_D^{h}v,v)=\int_\Omega \sum_{0\leq \alpha,\beta\leq 1}
g^{\alpha\beta}(X)\left( \nabla^h_\alpha-\frac14 ih
|g(X)|^{-1}\frac{\partial}{\partial X_\alpha}|g(X)|\right)
v(X)\times
\\ \times \overline{\left(\nabla^h_\beta -\frac14 ih
|g(X)|^{-1}\frac{\partial}{\partial X_\alpha}|g(X)|\right) v(X)}\,dX\\
-h\int_\Omega b(X)|v(X)|^2 dX +h^{1/2} \int_\Omega
(b(X)-b_0)|v(X)|^2\, dX\,.
\end{multline*}
Put
\[
q(v)=\int_\Omega \sum_{0\leq \alpha\leq 1} \left|\left(
\nabla^h_\alpha-\frac14 ih |g(X)|^{-1}\frac{\partial}{\partial
X_\alpha}|g(X)|\right) v(X)\right|^2 dX-hb_0\int_\Omega |v(X)|^2 \,dX\,.
\]
Then we have
\begin{multline}\label{e:qv}
|(\hat{P}_D^{h}v,v)-q(v)|\\ \leq \int_\Omega |t| \left|\left(
\nabla^h_0-\frac14 ih |g(X)|^{-1}\frac{\partial}{\partial
s}|g(X)|\right) v(X)\right|^2 \,dX\\
+C_1 h\int_\Omega t^2 |v(X)|^2 dX +C_2h^{1/2} \int_\Omega t^2
|v(X)|^2 \,dX\,.
\end{multline}
Consider
the Dirichlet realization $P_{mod,D}^{h}$ of the operator
\[
P_{mod}^{h}=\left(i h \frac{\partial}{\partial s}- b_0t\right)^2-h^2
\frac{\partial^2}{\partial t^2}-hb_0+ h^{1/2}\,
\left(\frac12\beta_2(s) t^2\right)
\]
in the space $L^2(\Omega,dX)$. So its quadratic form is given by
\[
(P_{mod}^{h}v,v)= q^{mod}(v)+h^{1/2} \int_\Omega\left(
\frac12\beta_2(s) t^2\right) |v(X)|^2 \,dX\,,
\]
where
\[
q^{mod}(v)=\int_\Omega \left|\left(i h \frac{\partial}{\partial
s}-b_0t\right)v(X)\right|^2 \, dX + h^2 \int_\Omega
\left|\frac{\partial}{\partial t}v(X)\right|^2 dX-hb_0\int_\Omega
|v(X)|^2 \,dX\,.
\]
So we have
\begin{multline}\label{e:P-P}
|(\hat{P}_D^{h}v,v)-(P_{mod}^{h}v,v)| \leq |q(v)-q^{mod}(v)|
\\ + \sum_\alpha \int_\Omega t^2 |\nabla^h_\alpha v(X)|^2 \,dX +
h\int_\Omega t^2 |v(X)|^2 \, dX +h^{1/2} \int_\Omega |t|^3|v(X)|^2\, dX\,.
\end{multline}
Finally, we have
\begin{multline}\label{e:q-q}
|q(v)-q^{mod}(v)|\\ \leq C(q(v))^{1/2}\left[h\left(\int_\Omega t^2
|v(X)|^2 dX\right)^{1/2}+\left(\int_\Omega t^6 |v(X)|^2
dX\right)^{1/2}\right].
\end{multline}

By Lemma~\ref{l:weight}, for any real $k\geq 0$, there exists
$C_k>0$ such that
\begin{equation}\label{e:Xu}
\| |t|^ku_h\|_{L^2(\Omega,dx_g)}\leq C_{k}(h^{k/2}+h^{(3k+1)/8}) \|
u_h\|_{L^2(\Omega,dx_g)},
\end{equation}
and, for any $\alpha =1,2$,
\begin{equation}\label{e:Xnablau}
\| |t|^k\nabla^h_\alpha u_h\|_{L^2(\Omega,dx_g)}\leq
C_{k}(h^{(k+1)/2}+h^{(3k+5)/8})\|u_h\|_{L^2(\Omega,dx_g)}.
\end{equation}

Put $v_h=|g(X)|^{1/4}u_h$. By \eqref{e:Xu} and \eqref{e:Xnablau},
for any real $k\geq 0$, there exists $C_k>0$ such that
\begin{equation}\label{e:Xv}
\| |t|^kv_h\|_{L^2(\Omega,dX)}\leq C_{k}(h^{k/2}+h^{(3k+1)/8}) \|
v_h\|_{L^2(\Omega,dX)},
\end{equation}
and, for any $\alpha =1,2$,
\begin{multline}\label{e:Xnablav}
\left( \int_\Omega |t|^{2k} \left|\left(i h \frac{\partial}{\partial
s}-b_0t\right)v_h(X)\right|^2 dX\right)^{1/2}+ \left(\int_\Omega
|t|^{2k} \left|h \frac{\partial}{\partial t}v_h(X)\right|^2
dX\right)^{1/2}\\ \leq
C_{k}(h^{(k+1)/2}+h^{(3k+5)/8})\|v_h\|_{L^2(\Omega,dX)}.
\end{multline}

The estimates \eqref{e:Xv} and \eqref{e:Xnablav} allow us to show,
first, using \eqref{e:Vhj}, \eqref{e:upperPh} and \eqref{e:qv},
that,
\[
q(v_h)\leq Ch^{3/2}\|v_h\|_{L^2(\Omega,dX)}^2,
\]
then, using \eqref{e:q-q}, that
\[
|q(v_h)-q^{mod}(v_h)|\leq Ch^2\|v_h\|_{L^2(\Omega,dX)}^2,
\]
and finally, using \eqref{e:P-P}, that
\begin{equation}\label{e:flatRn}
(P_{mod}^{h}v_h,v_h)\leq
(\lambda_0(P^h_D)+Ch^{15/8})\|v_h\|_{L^2(\Omega,dX)}^2.
\end{equation}

Let $\chi$ be a function from $C^\infty_c(\Omega)$ such that
$\chi=\chi(t)$ and $\chi\equiv 1$ in a neighborhood of $\gamma$. By
\eqref{e:Xv} and \eqref{e:Xnablav}, it follows that, for any
$k\in\NN$ there exists $C_k>0$ such that
\begin{equation}\label{e:chiv}
\|(1-\chi) v_h\|_{L^2(\Omega,dX)}+ \|\frac{\partial \chi}{\partial
t} v_h\|_{L^2(\Omega,dX)} \leq C_{k}h^{k}\| v_h\|_{L^2(\Omega,dX)}\,.
\end{equation}
and
\begin{multline}\label{e:chinablav}
\left( \int_\Omega (1-\chi(t)) \left|\left(i h
\frac{\partial}{\partial s}-b_0t\right)v_h(X)\right|^2
dX\right)^{1/2} \\ + \left(\int_\Omega (1-\chi(t)) \left| h
\frac{\partial}{\partial t}v_h(X)\right|^2 dX\right)^{1/2} \leq
C_{k}h^k\|v_h\|_{L^2(\Omega,dX)}\,.
\end{multline}

Using \eqref{e:chiv} and \eqref{e:chinablav}, it is easy to check
that, for any $k>0$, there exists $C_{k}>0$ such that
\begin{equation}\label{e:flat-chi}
|(P_{mod}^{h}(\chi v_h), \chi v_h)-(P_{mod}^{h}v_h, v_h)| \leq C_k
h^k \|v_h\|^2\,.
\end{equation}

Consider the self-adjoint realization of the operator $P_{mod}^{h}$
in $L^2(\RR^2,dX)$. We will keep  the same notation $P_{mod}^{h}$
for this operator. Put $w_h=\chi |g(X)|^{1/4}u_h$. By
\eqref{e:flatRn} and \eqref{e:flat-chi}, it follows that
\[
(P_{mod}^{h}w_h,w_h)\leq
(\lambda_0(P^h_D)+C_0h^{15/8})\|w_h\|_{L^2(\Omega,dX)}^2\,.
\]
This immediately implies that
\begin{equation}\label{e:flat-D-Rn}
\lambda_0(P_{mod}^{h}) \leq \lambda_0(P^h_D)+C_0h^{15/8}\,.
\end{equation}
Consider the operator
\[
P_{0}^{h}=\left(i h \frac{\partial}{\partial s}- b_0t\right)^2-h^2
\frac{\partial^2}{\partial t^2}-hb_0+ h^{1/2}\, \left(\frac12\mu_0
t^2\right)\,.
\]
Recall that the eigenvalues of the Schr\"odinger operator with
constant magnetic field and positive quadratic potential in $\RR^n$
can be computed explicitly. More precisely (see for instance
\cite[Theorem 2.2]{MatUeki}), the eigenvalues of the operator
\[
H_{b,K}=\left(i \frac{\partial}{\partial x}-\frac12
by\right)^2+\left(i\frac{\partial}{\partial y}+ \frac12
bx\right)^2+K_{11}x^2+2K_{12}xy+K_{22}y^2
\]
are given by
\[
\lambda_{n_1n_2}=(2n_1+1)s_1+(2n_2+1)s_2,\quad n_1,n_2\in \NN\,,
\]
where
\[
s_1=\frac{1}{\sqrt{2}}\left[t_K+b^2-[(t_K+b^2)^2-4d_K]^{1/2}\right]^{1/2}\,,
\]
\[
s_2=\frac{1}{\sqrt{2}}\left[t_K+b^2+[(t_K+b^2)^2-4d_K]^{1/2}\right]^{1/2}\,,
\]
and
\[
t_K=\operatorname{Tr}K=K_{11}+K_{22}\, , \quad d_K=\det
K=K_{11}K_{22}-K_{12}^2\,.
\]
Applying this formula to the operator $P_{0}^{h}$, we obtain that
its eigenvalues have the form ($t_K=\frac12\mu_0$, $d_K=0$):
\[
\lambda_{n_1n_2}=(2n_1+1)s_1+(2n_2+1)s_2-hb_0\,,\quad n_1,n_2\in \NN\,,
\]
where
\[
s_1 = 0\,,
\]
and
\begin{align*}
s_2 = h\left[h^{1/2} t_K+b_0^2\right]^{1/2} =  hb_0+\frac{t_K}{2b_0}  h^{3/2}+\mathcal O(h^{2})\, .
\end{align*}

For the lowest eigenvalue $\lambda_0(P_{0}^{h})$ of $P_{0}^{h}$, we
obtain
\begin{equation}\label{e:flat}
\lambda_0(P_{0}^{h})= \frac{\mu_0}{4b_0}h^{3/2}+\mathcal O(h^{2})\,
.
\end{equation}
Observing that $\lambda_0(P_{mod}^{h})\geq \lambda_0(P_{0}^{h})$ and
combining \eqref{e:flat-D-Rn} and \eqref{e:flat}, we immediately
complete the proof of Proposition~\ref{t:PAV}.

\section{Miniwells}\label{s:miniwells}
\subsection{Main statement}~\\
Fix some $k\in\NN$ and additionally assume that there exists a
unique minimum $x_0\in \gamma$ of the function
\begin{equation}\label{e:defVk}
V_k(x):=(2k^2+2k+1)\frac{\beta_2(x)}{4b_0}+\frac12
\left(k^2+k\right) R(x)
\end{equation}
on $\gamma$, which is nondegenerate, that is satisfying, for all $x\in \gamma$ in some
neighborhood of $x_0$,
\begin{equation}\label{e:wells1}
Cd(x,x_0)^2< V_k(x) - V_k(x_0)< C^{-1}d(x,x_0)^2\,.
\end{equation}

The purpose of this section is to give the following more precise
construction of approximate eigenvalues of the operator $H^h$.

\begin{thm}\label{t:qmodes0}
Under current assumptions, for any natural $j$, there exist
$u^h_{jk}\in C^\infty_c(M)$, $C_{jk}>0$ and $h_{jk}>0$ such that
\[
(u^h_{j_1k},u^h_{j_2k}) =\delta_{j_1j_2}+ \mathcal O_{j_1,j_2,k}(h)
\]
and, for any $h\in ]0,h_{jk}]$,
\[
\|H^hu^h_{jk}- \mu_{jk}^h u^h_{jk}\|\leq
C_{jk}h^{11/4}\|u^h_{jk}\|,\,
\]
where
\begin{equation}\label{e:mujk}
\mu_{jk}^h= \mu_{j,k,0}h+\mu_{j,k,4} h^2+ \mu_{j,k,6} h^{5/2},
\end{equation}
with
\[
\mu_{j,k,0} = (2k+1)b_0,\quad \mu_{j,k,4} = V_k(x_0),
\]
and
\[
\mu_{j,k,6}=\frac{1}{2b_0}
V^{\prime\prime}_k(x_0)^{1/2}\beta_2(x_0)^{1/2}(2k+1)^{1/2}(2j+1).
\]
\end{thm}

Here and below the derivative means the derivative with respect to
the natural parameter on $\gamma$.

As above, denote
\[
\mu_0:=\inf_{x\in\gamma} \beta_2(x).
\]
When $k=0$, we get:
\begin{cor}~\\ \label{c:qmodes0}
Assume that there exists a unique minimum point $x_0\in \gamma$ of
the function $\beta_2$ on $\gamma$, which is nondegenerate~:
\[
\mu_2:= \beta^{\prime\prime}_2(x_0)>0.
\]
 For any natural $j$, there exist $u^h_{j}\in C^\infty_c(M)$,
$C_{j}>0$ and $h_{j}>0$ such that
\[
(u^h_{j_1},u^h_{j_2}) =\delta_{j_1j_2}+ \mathcal O_{j_1,j_2}(h)
\]
and, for any $h\in ]0,h_{j}]$,
\[
\|H^hu^h_{j}- \mu_{j}^h u^h_{j}\|\leq C_{j}h^{11/4}\|u^h_{j}\|,\,
\]
where
\[
\mu_{j}^h= hb_0 + h^2\,\frac{\mu_0}{4b_0}+
h^{5/2}\,\frac{(\mu_0\mu_2)^{1/2}}{4b_0^{3/2}}(2j+1).
\]
\end{cor}

\begin{cor}~\\ \label{c:upper-lambda1}
Under the assumptions of Corollary~\ref{t:qmodes0}, for any natural
$j$, there exist $C_j$ and $h_j>0$, such that for any $h\in ]0,h_j]$
\[
\lambda_j(H^h)\leq hb_0 + h^2\,\frac{\mu_0}{4b_0}+
h^{5/2}\,\frac{(\mu_0\mu_2)^{1/2}}{4b_0^{3/2}}(2j+1) + C_j
h^{11/4}\,.
\]
\end{cor}

\begin{rem}
We conjecture that  under the assumptions of
Corollary~\ref{t:qmodes0}
$$
\lambda_0(H^h)=  hb_0 + h^2\,\frac{\mu_0}{4b_0}+
h^{5/2}\,\frac{(\mu_0\mu_2)^{1/2}}{4b_0^{3/2}}+ o ( h^{5/2})\,.
$$
\end{rem}

\subsection{Expanding operators in fractional powers of
$h$}~\\
 The approximate eigenfunctions $\Phi^h_k\in C^\infty(M)$, which
we are going to construct, will be supported in a small neighborhood
of $\gamma$. As in Section~\ref{s:expand}, we will use the normal
coordinate system $(s,t)$ in a tubular neighborhood $U$ of $\gamma$
with coordinates $X=(X_0,X_1)$ with $X_0=s$ and $X_1=t$. We make use
of notation of Section~\ref{s:expand} and need to consider further
terms in the asymptotic expansions of that section.

Thus, instead of \eqref{e:expG}, we write
\begin{equation}\label{e:expG-mod}
\begin{aligned}
g_{00}(s,t) & =1+a_{1}(s)t+{a}_{2}(s)t^2+{a}_{3}(s)t^3+\mathcal
O(t^4),\\ g_{11}(s,t) & = 1, \quad g_{01}(s,t) =0.
\end{aligned}
\end{equation}
Then \eqref{e:expG00} takes the form
\begin{multline}\label{e:expG00-mod}
g^{00}(s,t)=1-{a}_{1}(s)t-({a}_{2}(s)-{a}_{1}(s)^2)t^2\\
-({a}_{3}(s)-2{a}_{1}(s)a_2(s)+a_1(s)^3)t^3 +\mathcal O(t^4)\,, \quad
t\to 0\,.
\end{multline}

For the components of the magnetic potential, we assume that
\[
A_1(s,t)\equiv 0\,,
\]
and, instead of \eqref{e:A0}, we get
\begin{multline}\label{e:A0-mod}
A_0(s,t)=\alpha_1-b_0t-\frac14{a}_{1}(s)b_0t^2\\
-\frac16\left(\beta_2(s)+b_0({a}_{2}(s)-\frac14
{a}_{1}(s)^2)\right)t^3+\mathcal O(t^4)\,, \quad t\to 0\,,
\end{multline}
where $\alpha_1$ is some constant. Without loss of generality, we
can locally  (after a gauge transformation) take $\alpha_1=0$.

We have
\begin{equation}\label{e:Gamma-mod}
\begin{aligned}
\Gamma^0 =& -\frac12a^\prime_{1}(s)t+\left(-\frac12a^\prime_{2}(s) +
{a}_{1}(s)a^\prime_1(s)\right)t^2+\mathcal O(t^3),\\
\Gamma^1 = & \frac12{a}_{1}(s)+\left({a}_{2}(s)-\frac12
{a}_{1}(s)^2\right)t \\ & + \left(\frac32{a}_{3}(s)-\frac32
a_1(s)a_2(s)+\frac{1}{2}a_1(s)^3\right)t^2 +\mathcal O(t^3)\,.
\end{aligned}
\end{equation}

We move the operators into the Hilbert space $L^2(\RR^n)$ equipped
with the Euclidean inner product. For the operators
$\hat{\nabla}^h_0$ and $\hat{\nabla}^h_1$ defined in
\eqref{e:defNabla}, we obtain the following expansions:
\begin{equation}\label{e:Nabla-mod}
\begin{aligned}
\hat{\nabla}^h_0 = & {\nabla}^h_0-\frac i4 h a^\prime_1(s)t-
\frac i4 h \left(a^\prime_2(s)- a_1(s)a^\prime_1(s)\right)t^2+\mathcal O(t^3)\,,\\
\hat{\nabla}^h_1
= & -hD_t+ih\Big[-\frac14 a_1(s)-\frac12(a_2(s)-\frac12 a_1(s)^2)t\\
& -\frac12\left(\frac32{a}_{3}(s)-\frac32
a_1(s)a_2(s)+\frac{1}{2}a_1(s)^3\right)t^2+\mathcal O(t^3)\Big]\,.
\end{aligned}
\end{equation}

We will assume that the minimum point $x_0\in \gamma$ corresponds to
$s=0$. Then the condition \eqref{e:wells1} implies that
\[
\frac{\beta^{\prime\prime}_2(0)}{4b_0}(2k^2+2k+1)-
\left(k^2+k\right)\left(a^{\prime\prime}_2(0)-\frac12a_1(0)a^{\prime\prime}_1(0)
-\frac12(a^\prime_1(0))^2\right)>0\,.
\]
We also have
\begin{equation}\label{e:wells1-local}
\frac{\beta^\prime_2(0)}{4b_0}(2k^2+2k+1)-
\left(k^2+k\right)\left(a^\prime_2(0)-\frac12a_1(0)a^\prime_1(0)\right)=0\,.
\end{equation}

Now we use the scaling $t=h^{1/2}t_1$, $s=h^{1/4}s_1$. Similarly to
Section~\ref{s:expand}, we will only apply our operator on functions
which are a product of cut-off functions with functions of the form
of linear combinations of terms like $h^\nu w(h^{-1/4}s,h^{-1/2}t)$
with $w$ in $C^\infty(S^1)\otimes \cS(\RR_t)$ supported in a small
neighborhood of $(0,0)\in S^1\times \RR_t$. These functions are
consequently $\mathcal O(h^\infty)$ outside a fixed neighborhood of
$(0,0)$. We will start by doing the computations formally in the
sense that everything is determined modulo $\mathcal O(h^\infty)$,
and any smooth function will be replaced by its Taylor's expansion
at $(0,0)$.

From \eqref{e:expG00-mod}, \eqref{e:A0-mod}, \eqref{e:Gamma-mod},
\eqref{e:Nabla-mod}, we derive the following expansions:
\begin{equation}\label{e:nabla0-hat}
\begin{aligned}
\hat\nabla^h_0= & -b_0 h^{1/2}t_1 + i h^{3/4}
\frac{\partial}{\partial s_1} -\frac14a_1(0)b_0
ht^2_1-\frac14 a^\prime_1(0)b_0 h^{5/4}s_1t^2_1\\
& +h^{3/2}\left[-\frac16\left(\beta_2(0)+b_0(a_2(0)-\frac14
a_1(0)^2)\right)t_1^3-\frac i4 a^\prime_1(0)t_1-\frac18
a^{\prime\prime}_1(0)b_0 s^2_1t^2_1\right]\\ &
+h^{7/4}\Big[-\frac16\left[\beta^\prime_2(0)+b_0\left(a^\prime_2(0)
-\frac12a_1(0)a^\prime_1(0)\right)\right]s_1t_1^3\\ & -\frac i4
a^{\prime\prime}_1(0)s_1t_1-\frac{1}{24}
a^{\prime\prime\prime}_1(0)b_0
s^3_1t^2_1\Big]\\
&
+h^{2}\Big[-\frac{1}{12}\left(\beta^{\prime\prime}_2(0)+b_0(a^{\prime\prime}_2(0)-\frac12
a_1(0)a^{\prime\prime}_1(0)
-\frac12(a^\prime_1(0))^2)\right)s_1^2t_1^3
\\ & -\frac i8 a^{\prime\prime\prime}_1(0)s_1^2 t_1-\frac{1}{96} a^{IV}_1(0)b_0
s^4_1t^2_1- \frac i4 \left(a^\prime_2(0)- a_1(0)a^\prime_1(0)\right)
t_1^2\Big] +\mathcal O(h^{9/4})\,.
\end{aligned}
\end{equation}

\begin{equation}\label{e:nabla1-hat}
\begin{aligned}
\hat{\nabla}^h_1= &  -h^{1/2} D_{t_1}-\frac i4
 h a_1(0)-\frac i4 h^{5/4} a^\prime_1(0)s_1\\
& -h^{3/2}\left[\frac i2 (a_2(0)-\frac12a_1(0)^2)t_1+\frac i8
 a^{\prime\prime}_1(0)s^2_1\right]\\ &
-h^{7/4}\left[\frac12i(a^\prime_2(0)-a_1(0)a^\prime_1(0))s_1t_1+\frac{i}{24}
 a^{\prime\prime\prime}_1(0)s^3_1\right]\\ &
-h^{2}\Big[\frac i2 \left(\frac32{a}_{3}(0)-\frac32
a_1(0)a_2(0)+\frac{1}{2}a_1(0)^3\right)t_1^2\\ & + \frac i4
(a^{\prime\prime}_2(0)-a_1(0)a^{\prime\prime}_1(0)-(a^\prime_1(0))^2)s^2_1t_1+\frac{i}{96}
a^{IV}_1(0)s^4_1\Big]+\mathcal O(h^{9/4})\,.
\end{aligned}
\end{equation}

\begin{equation}\label{e:G00-hat}
\begin{aligned}
g^{00}=& 1-{a}_{1}(0)h^{1/2}t_1-a^\prime_{1}(0)s_1h^{3/4}t_1\\
&
+h\left[-({a}_{2}(0)-{a}_{1}(0)^2)t^2_1-\frac12a^{\prime\prime}_{1}(0)s^2_1t_1\right]\\
&+h^{5/4}\left[-\frac16a^{\prime\prime\prime}_{1}(0)s^3_1t_1-
(a^\prime_{2}(0)-2a_{1}(0)a^\prime_1(0))s_1t^2_1\right]\\
& +h^{3/2} \Big[-({a}_{3}(0)-2{a}_{1}(0)a_2(0)+a_1(0)^3)t^3_1 \\
&-\frac12(a^{\prime\prime}_{2}(0)-2a_{1}(0)a^{\prime\prime}_1(0)
-2a^\prime_1(0)^2)s^2_1t^2_1-\frac{1}{24}a^{IV}_{1}(0)s^4_1t_1\Big]+\mathcal
O(h^{7/4})\,.
\end{aligned}
\end{equation}

\begin{equation}\label{e:Gamma0-hat}
\begin{aligned}
\Gamma^0
=& h^{1/2}\left[-\frac12a^\prime_{1}(0)t_1\right]+h^{3/4}\left[-\frac12a^{\prime\prime}_{1}(0)s_1t_1\right]\\
& +h\left[\left(-\frac12a^\prime_{2}(0) +
{a}_{1}(0)a^\prime_1(0)\right)t_1^2-\frac14a^{\prime\prime\prime}_{1}(0)s^2_1t_1\right]+\mathcal
O(h^{5/4})\,.
\end{aligned}
\end{equation}

\begin{equation}\label{e:Gamma1-hat}
\begin{aligned}
\Gamma^1 = &
\frac12{a}_{1}(0)+h^{1/4}\frac12a^\prime_{1}(0)s_1+h^{1/2}\left[\left({a}_{2}(0)-\frac12
{a}_{1}(0)^2\right)t_1+\frac14a^{\prime\prime}_{1}(0)s_1^2\right] \\
& +h^{3/4}\left[\frac{1}{12}a^{\prime\prime\prime}_{1}(0)s_1^3+
\left(a^\prime_{2}(0)-
{a}_{1}(0)a_1^\prime(0)\right)s_1t_1\right]\\
& + h\Big[ \left(\frac32{a}_{3}(0)-\frac32
a_1(0)a_2(0)+\frac{1}{2}a_1(0)^3\right)t_1^2\\ & +
\frac12\left(a^{\prime\prime}_{2}(0)-
{a}_{1}(0)a_1^{\prime\prime}(0)-(a^\prime_1(0))^2\right)s^2_1t_1+\frac{1}{48}a^{IV}_{1}(0)s_1^4\Big]
+\mathcal O(h^{5/4})\,.
\end{aligned}
\end{equation}
Next, we make the partial Fourier transform $F_{s_1\to \sigma_1}$ in
the $s_1$-variable and the translation
\[
t_2=t_1+h^{1/4}\frac{\sigma_1}{b_0}, \quad \sigma_2=\sigma_1\,,
\]
and expand the operator $\hat{H}^h$ defined by \eqref{e:def-hatH} in
powers of $h^{1/4}$. Using \eqref{e:form-hatH},
\eqref{e:nabla0-hat}, \eqref{e:nabla1-hat}, \eqref{e:G00-hat},
\eqref{e:Gamma0-hat}, \eqref{e:Gamma1-hat}, after routine
computations, we obtain
\[
\hat{H}^h=hH^{h,0}\,,
\]
where
\[
H^{h,0}=P_0+h^{1/4}P_1+h^{1/2}P_2+h^{3/4}P_3+hP_4+h^{5/4}P_5+h^{3/2}P_6+\mathcal O(h^{7/4})\,,
\]
with
\begin{align*}
P_0=& - \frac{{\partial}^2}{\partial t_2^2} + b^2_0t^2_2\,,\\
P_1=& 0\,, \\
P_2=&-\frac12 a_1(0)b^2_0t_2^3\,,\\
P_3=& \frac12a^\prime_1(0)b_0^2 t_2^3D_{\sigma_2} \,,\\
P_4 =& \frac13\left(\beta_2(0)-2b_0\left(a_2(0)-\frac{23}{32}
a_1(0)^2\right)\right)b_0t^4_2+\frac 12
\left(a_2(0)-\frac38a_1(0)^2\right)\\
& +\frac12 a_1(0)t_2\sigma_2^2 + \frac12a^\prime_1(0)b_0
\left(t^3_2D_{t_2}-\frac{3}{2}it^2_2\right)-\frac14 a^{\prime\prime}_1(0)b^2_0 t_2^3 D_{\sigma_2}^2 \,,\\
P_5
=& -\left(\beta_2(0)-b_0(a_2(0)-\frac12 a_1(0)^2)\right)t^3_2\sigma_2\\
&-\frac13
b_0\left[\beta^\prime_2(0)-2b_0\left(a^\prime_2(0)-\frac{23}{16}a_1(0)a^\prime_1(0)\right)\right]t_2^4D_{\sigma_2}+
\frac{1}{12}a^{\prime\prime\prime}_1(0)b^2_0t_2^3 D_{\sigma_2}^3 \\
& -\frac14 a^{\prime\prime}_1(0)
b_0\left(2t_2^3D_{t_2}D_{\sigma_2}-3it_2^2D_{\sigma_2}\right) -
\frac12a^\prime_1(0)
t_2\sigma_2(\sigma_2D_{\sigma_2}-2t_2D_{t_2})\\
& -\frac12
\left(a^\prime_2(0)-\frac34a_1(0)a^\prime_1(0)\right)D_{\sigma_2}\,.
\end{align*}
The formula for $P_6$ is quite long, but we will only need the even
part of $P_6$ given by
\begin{align*}
P^+_6 =&
\Big[\left(\frac{\beta^{\prime\prime}_2(0)}{6b_0}-\frac{1}{3}\left(a^{\prime\prime}_2(0)
-\frac{23}{16}a_1(0)a^{\prime\prime}_1(0)
-\frac{23}{16}(a^\prime_1(0))^2\right)\right)b^2_0t_2^4
\\
& + \left(\frac 14 a^{\prime\prime}_2(0)-\frac
{3}{16}a_1(0)a^{\prime\prime}_1(0)-\frac
{3}{16}(a^\prime_1(0))^2\right)\\ &
-\frac{1}{4}a^{\prime\prime\prime}_1(0)b_0\left(t_2^3D_{t_2}-\frac32it_2^2\right)\Big] D_{\sigma_2}^2\\
& +\left[\frac{\beta_2(0)}{b_0}t^2_2-\frac38 a_1(0)^2t^2_2 -
\frac{1}{4b_0} a^\prime_1(0) (2t_2D_{t_2}-i)\right]\sigma_2^2\,.
\end{align*}

\subsection{Construction of approximate eigenfunctions}~\\
Now we look for an approximate eigenfunction $\varphi^h$ of $H^{h,0}$
and the corresponding approximate eigenvalue $\lambda^h$ as formal
expansions in powers of $h^{1/4}$:
\[
\varphi^h(t_2,\sigma_2)\sim
\sum_{\ell=0}^{\infty}h^{\ell/4}\varphi_\ell (t_2,\sigma_2)\,, \quad
\lambda^h\sim \sum_{\ell=0}^{\infty}h^{\ell/4}\lambda_\ell\,.
\]
In the next steps we will express the cancellation of the
coefficients of $h^{{\ell}/{4}}$ in the formal expansion of
$(H^{h,0}-\lambda^h)\varphi^h$
for $\ell= 0,1,\ldots,6$\,.\\
{\bf Step 1}~\\
At the first step we get
\[
P_0\varphi_0=\lambda_0\varphi_0\,,
\]
so we take
\[
\lambda_0=(2k+1)b_0\,,\quad
\varphi_0(t_2,\sigma_2)=\chi_0(\sigma_2)\psi_k(t_2)\,,
\]
with the given $k$.
\medskip \par
\noindent {\bf Step 2}~\\
The second cancellation gives:
\[
(P_0-\lambda_0)\varphi_1=\lambda_1 \varphi_0 - P_1 \varphi_0=\lambda_1\varphi_0\,,
\]
so we take
\[
\lambda_1=0\,,\quad \varphi_1=0\,.
\]
{\bf Step 3}~\\
The next equation reads:
\[
(P_0-\lambda_0)\varphi_2=\lambda_2\varphi_0-P_2\varphi_0-P_1\varphi_1=\lambda_2\varphi_0-P_2\varphi_0\,.
\]
The computation of $P_2\varphi_0$ gives:
\begin{align*}
P_2\varphi_0=& -\frac12
a_1(0)b^2_0\chi_0(\sigma_2)t_2^3\psi_k(t_2)\\
= & -\frac{1}{16}
a_1(0)b^{1/2}_0\chi_0(\sigma_2)(\psi_{k+3}(t_2)+(6k+6)\psi_{k+1}(t_2)\\
& +12k^2\psi_{k-1}(t_2)+8k(k-1)(k-2)\psi_{k-3}(t_2))\,.
\end{align*}
Thus, we have
\[
\lambda_2=0\,,
\]
and
\begin{align*}
(P_0-\lambda_0)\varphi_2= & \frac{1}{16}
a_1(0)b^{1/2}_0\chi_0(\sigma_2)(\psi_{k+3}(t_2)+(6k+6)\psi_{k+1}(t_2)\\
& +12k^2\psi_{k-1}(t_2)+8k(k-1)(k-2)\psi_{k-3}(t_2))\,.
\end{align*}
We can take
\begin{align*}
\varphi_2 = & \frac{1}{16}
a_1(0)b^{-1/2}_0\chi_0(\sigma_2)(\frac{1}{6}\psi_{k+3}(t_2)+(3k+3)\psi_{k+1}(t_2)\\
& -6k^2\psi_{k-1}(t_2)-\frac{4}{3}k(k-1)(k-2)\psi_{k-3}(t_2))\,.
\end{align*}
{\bf Step 4}~\\
As fourth equation we obtain:
\[
(P_0-\lambda_0)\varphi_3=\lambda_3\varphi_0-P_3\varphi_0\,.
\]
The computation of $P_3\varphi_0$ gives:
\begin{align*}
P_3\varphi_0=& - \frac12ia^\prime_1(0)b_0^2\chi^\prime_0(\sigma_2)t_2^3\psi_k(t_2)\\
= & - \frac{1}{16}i
a^\prime_1(0)b^{1/2}_0\chi^\prime_0(\sigma_2)(\psi_{k+3}(t_2)+(6k+6)\psi_{k+1}(t_2)\\
& +12k^2\psi_{k-1}(t_2)+8k(k-1)(k-2)\psi_{k-3}(t_2))\,.
\end{align*}
Thus, we have
\[
\lambda_3=0\,,
\]
and
\begin{align*}
(P_0-\lambda_0)\varphi_3=  & \frac{1}{16}i
a^\prime_1(0)b^{1/2}_0\chi^\prime_0(\sigma_2)(\psi_{k+3}(t_2)+(6k+6)\psi_{k+1}(t_2)\\
& +12k^2\psi_{k-1}(t_2)+8k(k-1)(k-2)\psi_{k-3}(t_2))\,.
\end{align*}
We can take
\begin{align*}
\varphi_3 = & \frac{1}{16}i
a^\prime_1(0)b^{-1/2}_0\chi^\prime_0(\sigma_2)(\frac{1}{6}\psi_{k+3}(t_2)+(3k+3)\psi_{k+1}(t_2)\\
& -6k^2\psi_{k-1}(t_2)-\frac{4}{3}k(k-1)(k-2)\psi_{k-3}(t_2))\,.
\end{align*}
{\bf Step 5}~\\
The fifth equation reads:
\[
(P_0-\lambda_0) \varphi_4=\lambda_4 \varphi_0 - P_4\varphi_0- P_2 \varphi_2\,.
\]
The computation of $P_4\varphi_0$ gives:
\begin{align*}
P_4\varphi_0=&
\frac13\left(\beta_2(0)-2b_0\left(a_2(0)-\frac{23}{32}
a_1(0)^2\right)\right)b_0t^4_2\psi_k(t_2)\chi_0(\sigma_2)\\ & +\frac
12 \left(a_2(0)-\frac38a_1(0)^2\right)\psi_k(t_2)\chi_0(\sigma_2)
+\frac12 a_1(0)t_2\psi_k(t_2)\sigma_2^2\chi_0(\sigma_2)\\
& + \frac12a^\prime_1(0)b_0
\left(t^3_2D_{t_2}-\frac{3}{2}it^2_2\right)\psi_k(t_2)\chi_0(\sigma_2)+\frac14
a^{\prime\prime}_1(0)b^2_0
t_2^3\psi_k(t_2)\chi^{\prime\prime}_0(\sigma_2)\,.
\end{align*}
Multiplying by $\psi_k$ and integrating with respect to $t_2$ gives:
\begin{multline*}
\frac{1}{\|H_k\|^2}\int P_4\varphi_0 (\sigma_2,t_2)\psi_k(t_2)\,dt_2\\
\begin{aligned}
=&
\left(\frac{1}{4b_0}\beta_2(0)-\frac{1}{2}\left(a_2(0)-\frac{23}{32}
a_1(0)^2\right)\right)(2k^2+2k+1) \chi_0(\sigma_2)\\
& +\frac 12 \left(a_2(0)-\frac38a_1(0)^2\right)\chi_0(\sigma_2)\,.
\end{aligned}
\end{multline*}
Compute
\begin{align*}
P_2\varphi_2=& \frac{1}{32}
a_1(0)^2b^{3/2}_0\chi_0(\sigma_2)(\frac{1}{6}t_2^3\psi_{k+3}(t_2)+(3k+3)t_2^3\psi_{k+1}(t_2)\\
&
-6k^2t_2^3\psi_{k-1}(t_2)-\frac{4}{3}k(k-1)(k-2)t_2^3\psi_{k-3}(t_2))\,.
\end{align*}
Multiplying by $\psi_k$ and integrating with respect to $t_2$ gives:
\[
\frac{1}{\|H_k\|^2}\int P_2\varphi_2 (\sigma_2,t_2)\psi_k(t_2)\,dt_2 =
- \left(\frac{15}{32}k^2+\frac{15}{32}k+
\frac{11}{64}\right)a_1(0)^2\chi_0(\sigma_2)\,.
\]
We obtain
\begin{align*}
\lambda_4\chi_0(\sigma_2) & -\frac{1}{\|H_k\|^2}\int
(P_4\varphi_0+P_2\varphi_2) \psi_k\, dt_2\\
=  \lambda_4\chi_0(\sigma_2) & - \frac{1}{4b_0}\beta_2(0)(2k^2+2k+1)
\chi_0(\sigma_2)-(k^2+k)\left(a_2(0)-\frac{1}{4} a_1(0)^2\right)
\chi_0(\sigma_2)\,,
\end{align*}
which holds for any function $\chi_0$, if we put
\[
\lambda_4 = \frac{1}{4b_0}\beta_2(0)(2k^2+2k+1) +
(k^2+k)\left(a_2(0)-\frac{1}{4} a_1(0)^2\right)\,.
\]

To find $\varphi_4$, we observe that the right hand side has the
form
\begin{multline*}
\lambda_4\varphi_0- P_4 \varphi_0- P_2 \varphi_2\\
\begin{aligned}
= & A_0(t_2)\chi_0(\sigma_2)-\frac{1}{4b_0^{1/2}}
a_1(0)(\psi_{k+1}(t_2)+2k\psi_{k-1}(t_2))\sigma_2^2\chi_0(\sigma_2)\\
& -\frac{1}{32} a^{\prime\prime}_1(0)b^{3/2}_0
\Big(\psi_{k+3}(t_2)+(6k+6)\psi_{k+1}(t_2)\\ & +12k^2\psi_{k-1}(t_2)
+8k(k-1)(k-2)\psi_{k-3}(t_2)\Big) \chi^{\prime\prime}_0(\sigma_2),
\end{aligned}
\end{multline*}
where $A_0(t_2)$ has a form
\[
A_0(t_2)=\sum_{m=-3}^3\alpha_m\psi_{k+2m}(t_2), \quad \alpha_0=0\,.
\]
Thus, we obtain
\begin{align*}
\varphi_4 = & B_0(t_2)\chi_0(\sigma_2)-\frac{1}{4b_0^{3/2}}
a_1(0)\left(\frac{1}{2}\psi_{k+1}(t_2)-k\psi_{k-1}(t_2)\right)\sigma_2^2\chi_0(\sigma_2)\\
& -\frac{1}{32} a^{\prime\prime}_1(0)b^{1/2}_0
\Big(\frac{1}{6}\psi_{k+3}(t_2)+(3k+3)\psi_{k+1}(t_2)\\
& -6k^2\psi_{k-1}(t_2)-\frac{4}{3}k(k-1)(k-2)\psi_{k-3}(t_2)\Big)
\chi^{\prime\prime}_0(\sigma_2)\,,
\end{align*}
where
\[
B_0(t_2)=\sum_{m=-3}^3\frac{\alpha_m}{2mb_0}\psi_{k+2m}(t_2)\,.
\]

{\bf Step 6}~\\
The sixth equation reads:
\[
(P_0-\lambda_0)\varphi_5=\lambda_5\varphi_0-P_5\varphi_0-P_3\varphi_2-P_2\varphi_3\,.
\]

We don't need an explicit formula for $\varphi_5$ but only its existence. Therefore, we
only find $\lambda_5$ from an orthogonality condition. We have
\begin{align*}
P_5\varphi_0
= & -\left(\beta_2(0)-b_0(a_2(0)-\frac12 a_1(0)^2)\right)t^3_2\psi_k(t_2)\sigma_2\chi_0(\sigma_2)\\
&+\frac13i
b_0\left[\beta^\prime_2(0)-2b_0\left(a^\prime_2(0)-\frac{23}{16}a_1(0)a^\prime_1(0)\right)\right]
t_2^4\psi_k(t_2)\chi^\prime_0(\sigma_2)\\ & -
\frac{1}{12}ia^{\prime\prime\prime}_1(0)b^2_0t_2^3\psi_k(t_2)\chi^{\prime\prime\prime}_0(\sigma_2)
 +\frac14 a^{\prime\prime}_1(0)
b_0\left(2t_2^3D_{t_2}-3it_2^2\right)\psi_k(t_2)\chi^\prime_0(\sigma_2)\\
&  - \frac12a^\prime_1(0)
t_2\sigma_2(\sigma_2D_{\sigma_2}-2t_2D_{t_2})\psi_k(t_2)\chi_0(\sigma_2)\\
& +\frac12i
\left(a^\prime_2(0)-\frac14a_1(0)a^\prime_1(0)\right)\psi_k(t_2)\chi^\prime_0(\sigma_2).
\end{align*}
Multiplying by $\psi_k$ and integrating with respect to $t_2$ gives:
\begin{multline*}
\frac{1}{\|H_k\|^2}\int P_5\varphi_0 (\sigma_2,t_2)\psi_k(t_2)\,dt_2\\
\begin{aligned}
= &i
\left[\frac{\beta^\prime_2(0)}{4b_0}(2k^2+2k+1)-\left(a^\prime_2(0)
-\frac{23}{16}a_1(0)a^\prime_1(0)\right)(k^2+k)\right]
\chi^\prime_0(\sigma_2) \\
& +\frac{11}{32}a_1(0)a^\prime_1(0)\chi^\prime_0(\sigma_2)\,.
\end{aligned}
\end{multline*}
Next,
\begin{align*}
P_3\varphi_2 = & -\frac{1}{32}i
a_1(0)a^\prime_1(0)b^{3/2}_0(\frac{1}{6}t_2^3\psi_{k+3}(t_2)+(3k+3)t_2^3\psi_{k+1}(t_2)\\
&
-6k^2t_2^3\psi_{k-1}(t_2)-\frac{4}{3}k(k-1)(k-2)t_2^3\psi_{k-3}(t_2))\chi^\prime_0(\sigma_2)\,.
\end{align*}

Multiplying by $\psi_k$ and integrating with respect to $t_2$ gives:
\[
\frac{1}{\|H_k\|^2}\int P_3\varphi_2 (\sigma_2,t_2)\psi_k(t_2)\,dt_2 =
-\left(\frac{15}{32}k^2+\frac{15}{32}k+\frac{11}{64}\right)i
a_1(0)a^\prime_1(0) \chi^\prime_0(\sigma_2)\,.
\]
Finally,
\begin{align*}
P_2\varphi_3 = & -\frac{1}{32}i
a_1(0)a^\prime_1(0)b^{3/2}_0(\frac{1}{6}t_2^3\psi_{k+3}(t_2)+(3k+3)t_2^3\psi_{k+1}(t_2)\\
&
-6k^2t_2^3\psi_{k-1}(t_2)-\frac{4}{3}k(k-1)(k-2)t_2^3\psi_{k-3}(t_2))\chi^\prime_0(\sigma_2)\,.
\end{align*}
Multiplying by $\psi_k$ and integrating with respect to $t_2$ gives:
\[
\frac{1}{\|H_k\|^2}\int P_2\varphi_3 (\sigma_2,t_2)\psi_k(t_2)\, dt_2 =
-\left(\frac{15}{32}k^2+\frac{15}{32}k+\frac{11}{64}\right)i
a_1(0)a^\prime_1(0) \chi^\prime_0(\sigma_2)\,.
\]
Using \eqref{e:wells1-local}, we obtain
\begin{multline*}
\frac{1}{\|H_k\|^2}\int
(P_5\varphi_0+P_3\varphi_2+ P_2\varphi_3)\psi_k\, dt_2\\
= i
\left[\frac{\beta^\prime_2(0)}{4b_0}(2k^2+2k+1)-\left(a^\prime_2(0)
-\frac{1}{2}a_1(0)a^\prime_1(0)\right)(k^2+k)\right]
\chi^\prime_0(\sigma_2) =  0\,.
\end{multline*}
Therefore, the orthogonality condition holds for any function
$\chi_0$, if we put
\[
\lambda_5=0\,.
\]

{\bf Step 7}~\\
The seventh equation reads:
\[
(P_0-\lambda_0)\varphi_6=\lambda_6\varphi_0-P_6\varphi_0-P_4\varphi_2-P_3\varphi_3-P_2\varphi_4\, .
\]
We have
\begin{align*}
P_6\varphi_0 = &
-\Big[\left(\frac{\beta^{\prime\prime}_2(0)}{6b_0}-\frac{1}{3}\left(a^{\prime\prime}_2(0)
-\frac{23}{16}a_1(0)a^{\prime\prime}_1(0)
-\frac{23}{16}(a^\prime_1(0))^2\right)\right)b^2_0t_2^4\psi_k(t_2)
\\
& + \left(\frac 14 a^{\prime\prime}_2(0)-\frac
{3}{16}a_1(0)a^{\prime\prime}_1(0)-\frac
{3}{16}(a^\prime_1(0))^2\right)\psi_k(t_2)\\ &
-\frac{1}{4}a^{\prime\prime\prime}_1(0)b_0\left(t_2^3D_{t_2}-\frac32it_2^2\right)\psi_k(t_2)\Big] \chi^{\prime\prime}_0(\sigma_2)\\
& +\left[\left(\frac{\beta_2(0)}{b_0}-\frac38
a_1(0)^2\right)t^2_2\psi_k(t_2) - \frac{1}{4b_0} a^\prime_1(0)
(2t_2D_{t_2}-i)\psi_k(t_2)\right]\sigma_2^2\chi_0(\sigma_2)\,.
\end{align*}
Multiplying by $\psi_k$ and integrating with respect to $t_2$ gives:
\begin{multline*}
\frac{1}{\|H_k\|^2}\int P_6\varphi_0 (\sigma_2,t_2)\psi_k(t_2)\, dt_2\\
\begin{aligned}
= & -\Big[\frac{\beta^{\prime\prime}_2(0)}{8b_0}(2k^2+2k+1)
-\frac{1}{2}(k^2+k)a^{\prime\prime}_2(0)\\
&
+\left(\frac{23}{32}k^2+\frac{23}{32}k+\frac{11}{64}\right)a_1(0)a^{\prime\prime}_1(0)
+\left(\frac{23}{32}k^2+\frac{23}{32}k+\frac{11}{64}\right)(a^\prime_1(0))^2
\Big] \chi^{\prime\prime}_0(\sigma_2)\\
& + \left(\frac{\beta_2(0)}{2b^2_0}(2k+1)-\frac{3}{16b_0}(2k+1)
a_1(0)^2\right)\sigma_2^2\chi_0(\sigma_2)\,.
\end{aligned}
\end{multline*}
Next, we have
\begin{align*}
P_2\varphi_4 = & -\frac12
a_1(0)b^2_0t_2^3\Big[B_0(t_2)\chi_0(\sigma_2)\\
& -\frac{1}{4b_0^{3/2}}
a_1(0)\left(\frac{1}{2}\psi_{k+1}(t_2)-k\psi_{k-1}(t_2)\right)\sigma_2^2\chi_0(\sigma_2)\\
& -\frac{1}{32} a^{\prime\prime}_1(0)b^{1/2}_0
\Big(\frac{1}{6}\psi_{k+3}(t_2)+(3k+3)\psi_{k+1}(t_2)
-6k^2\psi_{k-1}(t_2)\\ & -\frac{4}{3}k(k-1)(k-2)\psi_{k-3}(t_2)\Big)
\chi^{\prime\prime}_0(\sigma_2)\Big]\,.
\end{align*}
Observe that
\[
\langle t_2^3
B_0,\psi_k\rangle=\sum_{m=-3}^3\frac{\alpha_m}{2mb_0}\langle
t_2^3\psi_{k+2m},\psi_k\rangle=0\,.
\]
Multiplying by $\psi_k$ and integrating with respect to $t_2$ gives:
\begin{multline*}
\frac{1}{\|H_k\|^2}\int P_2\varphi_4 (\sigma_2,t_2)\psi_k(t_2)\,dt_2\
= \frac{3}{32b_0}(2k+1) (a_1(0))^2\sigma_2^2\chi_0(\sigma_2)\\ +
a_1(0) a^{\prime\prime}_1(0)b_0
\Big(\frac{15}{64}k^2+\frac{15}{64}k+\frac{11}{128}\Big)
\chi^{\prime\prime}_0(\sigma_2)\,.
\end{multline*}

Next,
\begin{align*}
P_3\varphi_3 = & \frac{1}{32}
(a^\prime_1(0))^2b^{3/2}_0t_2^3(\frac{1}{6}\psi_{k+3}(t_2)+(3k+3)\psi_{k+1}(t_2)\\
&
-6k^2\psi_{k-1}(t_2)-\frac{4}{3}k(k-1)(k-2)\psi_{k-3}(t_2))\chi^{\prime\prime}_0(\sigma_2)\,.
\end{align*}

Multiplying by $\psi_k$ and integrating with respect to $t_2$ gives:
\[
\frac{1}{\|H_k\|^2}\int P_3\varphi_3 (\sigma_2,t_2)\psi_k(t_2)\, dt_2
=\left(\frac{15}{32}k^2+\frac{15}{32}k+\frac{11}{64}\right)
(a^\prime_1(0))^2\chi^{\prime\prime}_0(\sigma_2)\,.
\]

Finally, the operator $P_4$ has the form
\[
P_4 = A_4(t_2)+\frac12 a_1(0)t_2\sigma_2^2
 -\frac14 a^{\prime\prime}_1(0)b^2_0 t_2^3 D_{\sigma_2}^2\,.
\]
Therefore,
\begin{align*}
P_4\varphi_2 =& B_4(t_2)\chi_0(\sigma_2) + \frac{1}{32}
(a_1(0))^2b^{-1/2}_0t_2\Big(\frac{1}{6}\psi_{k+3}(t_2)+(3k+3)\psi_{k+1}(t_2)\\
& -6k^2\psi_{k-1}(t_2)-\frac{4}{3}k(k-1)(k-2)\psi_{k-3}(t_2)\Big)\sigma_2^2 \chi_0(\sigma_2)\\
& + \frac{1}{64}
a_1(0)a^{\prime\prime}_1(0)b^{3/2}_0t_2^3\Big(\frac{1}{6}\psi_{k+3}(t_2)+(3k+3)\psi_{k+1}(t_2)\\
&
-6k^2\psi_{k-1}(t_2)-\frac{4}{3}k(k-1)(k-2)\psi_{k-3}(t_2)\Big)\chi^{\prime\prime}_0(\sigma_2)\,.
\end{align*}

It is easy to see that the term $B_4(t_2)$ has the form
\[
B_4(t_2)=\sum_{m=-4}^3\gamma_m\psi_{k+2m+1}(t_2)\,.
\]
Therefore, we have
\[
\langle B_4,\psi_k\rangle=0\,.
\]

Multiplying by $\psi_k$ and integrating with respect to $t_2$ gives:
\begin{multline*}
\frac{1}{\|H_k\|^2}\int P_4\varphi_2 (\sigma_2,t_2)\psi_k(t_2)\, dt_2\\
=\frac{3}{32b_0}(2k+1) (a_1(0))^2\sigma_2^2 \chi_0(\sigma_2)+
a_1(0)a^{\prime\prime}_1(0)\Big(\frac{15}{64}k^2+\frac{15}{64}k+\frac{11}{128}\Big)\chi^{\prime\prime}_0(\sigma_2)\,.
\end{multline*}
The orthogonality condition gives
\begin{align*}
\lambda_6\chi_0(\sigma_2) & -\frac{1}{\|H_k\|^2}\int
(P_6\varphi_0+P_4\varphi_2+ P_3\varphi_3+ P_4\varphi_2)\psi_k\, dt_2\\
=  \lambda_6\chi_0(\sigma_2) & +\frac12
\Big[\frac{\beta^{\prime\prime}_2(0)}{4b_0}(2k^2+2k+1)\\ &
-(k^2+k)\left( a^{\prime\prime}_2(0)
-\frac{1}{2}a_1(0)a^{\prime\prime}_1(0)-
\frac{1}{2}(a^\prime_1(0))^2\right) \Big]
\chi^{\prime\prime}_0(\sigma_2)\\ & -
\frac{\beta_2(0)}{2b^2_0}(2k+1)\sigma_2^2\chi_0(\sigma_2)\,,
\end{align*}
which has a nontrivial solution $\chi_0\,$, if $\lambda_6$ is an
eigenvalue of the operator
\begin{multline*}
D_k=-\frac12 \Big[\frac{\beta^{\prime\prime}_2(0)}{4b_0}(2k^2+2k+1)
-(k^2+k)\left( a^{\prime\prime}_2(0)
-\frac{1}{2}a_1(0)a^{\prime\prime}_1(0)-
\frac{1}{2}(a^\prime_1(0))^2\right) \Big] \frac{d^2}{d\sigma^2_2}\\
+ \frac{\beta_2(0)}{2b^2_0}(2k+1)\sigma_2^2\,.
\end{multline*}

Thus, we can take
\begin{multline*}
\lambda_6 =\frac{1}{2b_0}
\Big[\frac{\beta^{\prime\prime}_2(0)}{4b_0}(2k^2+2k+1)\\
-(k^2+k)\left( a^{\prime\prime}_2(0)
-\frac{1}{2}a_1(0)a^{\prime\prime}_1(0)-
\frac{1}{2}(a^\prime_1(0))^2\right) \Big]^{1/2}
\beta_2(0)^{1/2}(2k+1)^{1/2}(2j+1), \\ j\in\NN\,,
\end{multline*}
and
\[
\chi_0(\sigma_2)=\Psi_{jk}(\sigma_2)\,,
\]
where $\Psi_{jk}$ is the normalized eigenfunction of the operator
$D_k$ associated with the eigenvalue $\lambda_6$.
\medskip
\par
Thus, for any $j\in\NN$, we have constructed an approximate
eigenfunction $\varphi^h_{jk}$ of the operator $H^{h,0}$ in the form
\[
\varphi^h_{jk}(t_2,\sigma_2)= \sum_{\ell=0}^{6}h^{j/4}\varphi_\ell
(t_2,\sigma_2)\,, \quad \varphi_\ell^{jk}\in {\cS(\RR^2)}\,,
\]
such that
\[
\varphi_0 (t_2,\sigma_2)
=\frac{1}{\|H_k\|}\psi_{k}(t_2)\Psi_{jk}(\sigma_2)\,.
\]
with the corresponding approximate eigenvalue
\[
\lambda^{jk}(h)= \sum_{\ell=0}^{6}h^{\ell/4}\lambda^{jk}_\ell\,.
\]
Then we have
\[
H^{h,0}\varphi^h_{jk}-\lambda^{jk}(h)\varphi^h_{jk}=\mathcal
O(h^{7/4})\,.
\]
Observe that $\lambda^{jk}(h)=h^{-1}\mu^{jk}(h)$, where
$\mu^{jk}(h)$ are given by \eqref{e:mujk}.

The constructed functions $\varphi^h_{jk}$ have sufficient decay
properties. Therefore, by changing back to the original coordinates
and multiplying by a fixed cut-off function, we obtain the desired
approximate eigenfunctions $u^h_{jk}$, completing the proof of
Theorem~\ref{t:qmodes0}.

\section{Periodic case and spectral gaps}\label{s:gaps}
In this section, we apply the previous results to the problem of
existence of gaps in the spectrum of a periodic magnetic
Schr\"odinger operator. For related results on spectral gaps for
periodic magnetic Schr\"odinger operators, see \cite{luminy} and
references therein.

Let $ M$ be a two-dimensional noncompact oriented manifold of
dimension $n\geq 2$ equipped with a properly discontinuous action of
a finitely generated, discrete group $\Gamma$ such that $M/\Gamma$
is compact. Suppose that $H^1(M, \RR) = 0$, i.e. any closed $1$-form
on $M$ is exact. Let $g$ be a $\Gamma$-invariant Riemannian metric
and $\bf B$ a real-valued $\Gamma$-invariant closed 2-form on $M$.
Assume that $\bf B$ is exact and choose a real-valued 1-form $\bf A$
on $M$ such that $d{\bf A} = \bf B$. Write ${\bf B}=b dx_g$, where
$b\in C^\infty(M)$ and $dx_g$ is the Riemannian volume form. Let
\[
b_0=\min_{x\in M}|b(x)|\,.
\]
Assume that $b_0>0$ and there exist a (connected) fundamental domain
$\cF$ and a constant $\epsilon_0>0$ such that
\[
 |b(x)| \geq b_0+\epsilon_0, \quad x\in \partial\cF\,.
\]

We will consider the magnetic Schr\"odinger operator $H^h$ as an
unbounded self-adjoint operator in the Hilbert space $L^2(M)$. Using
the results of \cite{diff2006}, one can immediately derive from
Theorem~\ref{t:qmodes} the following result on existence of gaps in
the spectrum of $H^h$ in the semiclassical limit.

\begin{thm}~\\
Assume the set $\{x\in M : |b(x)|=b_0\}$ contains a
smooth curve $\gamma$ such that in some neighborhood of $\gamma$ we
have
\[
C^{-1}d(x,\gamma)^2\leq |b(x)| -b_0 \leq C
d(x,\gamma)^2\,
\]
with a constant $C>0$. Then, for any natural $k$ and $N$, there
exists $h_{k,N}>0$ such that the spectrum of $H^h$ in the interval
\[
[(2k+1) hb_0+ h^2m_k,(2k+1) hb_0+ h^2M_k]\,,
\]
where $[m_k,M_k]$ is the range of the function $V_k$ on $\gamma$
defined by \eqref{e:defVk}, has at least $N$ gaps for any $h\in ]0,
h_{k,N}]\,$.
\end{thm}

Similarly, Theorem~\ref{t:qmodes0} implies the following result.

\begin{thm}~\\
Assume the set $\{x\in M : |b(x)|=b_0\}$ contains a
smooth curve $\gamma$ such that in some neighborhood of $\gamma$ we
have
\[
C^{-1}d(x,\gamma)^2\leq  |b(x)| -b_0 \leq C
d(x,\gamma)^2\,
\]
with a constant $C>0$, and for some $k\in\NN$, there exists a unique
minimum $x_0\in \gamma$  of the function $V_k$ on $\gamma$ defined
by \eqref{e:defVk}, which is nondegenerate~:
\[
\delta_k=V^{\prime\prime}_k(x_0)>0\,.
\]
Then, for any natural $N$, there exist
\[
C_{k,N}>\frac{1}{2b_0} \delta_k^{1/2}\mu_0^{1/2}(2k+1)^{1/2}(2N+3)
\]
and $h_{k,N}>0$ such that the spectrum of $H^h$ in the interval
\[
[(2k+1) hb_0+ h^2m_k+h^{5/2}\frac{1}{2b_0}
\delta_k^{1/2}\mu_0^{1/2}(2k+1)^{1/2},(2k+1) hb_0+
h^2m_k+h^{5/2}C_{k,N}]
\]
has at least $N$ gaps for any $h\in ] 0, h_{k,N}]$\,.
\end{thm}

\appendix
\section{Some facts from geometry}
Here we refer to \cite{Gray} for more material.
Let $M$ be a Riemannian manifold. Denote by $\nabla$ the Levi-Civita
connection associated with the Riemannian metric $g$. Recall that a
connection in the tangent bundle of $M$ is a map
\[
\nabla : C^\infty(M,TM)\times C^\infty(M,TM)\to C^\infty(M,TM),\quad
(X,Y)\mapsto \nabla_XY\,,
\]
which satisfies
\begin{enumerate}
\item for any $X\in C^\infty(M,TM)$ the map
$\nabla_X : C^\infty(M,TM)\to C^\infty(M,TM)$ is linear;
\item for any $X_1,X_2\in C^\infty(M,TM)$ and  $Y\in C^\infty(M,TM)$
  \[
  \nabla_{X_1+X_2}s=\nabla_{X_1}Y+\nabla_{X_2}Y\,;
\]
  \item for any $f\in C^\infty(M)$, $X\in C^\infty(M,TM)$ and $Y\in C^\infty(M,TM)\,$
  \[
  \nabla_{fX}Y=f\nabla_XY\,;
\]
  \item for any $f\in C^\infty(M)$, $X\in C^\infty(M,TM)$ and $Y\in C^\infty(M,TM)$
  \[
\nabla_X(fY)=f\nabla_XY+X(f)Y\,.
  \]
\end{enumerate}

The Levi-Civita connection is the unique connection $\nabla$ in
$TM$, which is compatible with the Riemannian metric (with the
corresponding scalar product denoted by $(\cdot,\cdot)$) in the
sense that, for any $X,Y,Z\in C^\infty(M,TM)$
\[
X(Y,Z)=(\nabla_XY,Z)+(Y,\nabla_XZ)\,,
\]
and torsion-free, in the sense that, for any $X,Y\in C^\infty(M,TM)\,$,
\[
\nabla_XY-\nabla_YX-[X,Y]=0\,.
\]
The curvature of the connection $\nabla$ is the operator
\[
R(X,Y): C^\infty(M,TM)\to C^\infty(M,TM)\,,
\]
associated with any pair $X,Y\in C^\infty(M,TM)$ by the formula
\[
R(X,Y)Z=\nabla_X\nabla_YZ-\nabla_Y\nabla_XZ-\nabla_{[X,Y]}Z
\]
for any $Z\in C^\infty(M,TM)\,$.

In a local coordinate system $(x^1,x^2,\ldots,x^n)$, one can write
\[
R\left(\frac{\partial }{\partial x^j}, \frac{\partial }{\partial
x^k}\right)\frac{\partial }{\partial x^\ell}=\sum_{i=1}^nR^i_{\ell
jk}\frac{\partial }{\partial x^i}\,,
\]
where
\[
R^i_{\ell jk}=\frac{\partial \Gamma^i_{k\ell}}{\partial
x^j}-\frac{\partial \Gamma^i_{j\ell}}{\partial x^k}+\sum_{m=1}^n
\Gamma^m_{k\ell}\Gamma^i_{jm}-\sum_{m=1}^n\Gamma^m_{j\ell}\Gamma^i_{km}\,.
\]
The Ricci tensor of $M$ is given by
\[
{\rm Ric}_{ab}=\sum_{i=1}^nR^i_{aib}\,.
\]
The scalar curvature of $M$ is defined by
\[
\kappa = \sum_{ab}g^{ab}{\rm Ric}_{ab\,}.
\]
The Riemannian curvature tensor is a $R$ of type (4,0), defined in
local coordinates as
\[
R_{ijk\ell}=\sum_{m=1}^ng_{im}R^m_{jk\ell}\,.
\]
One can also give its invariant definition. The corresponding map
\[
R : {\mathcal X}(M)\times {\mathcal X}(M)\times {\mathcal
X}(M)\times {\mathcal X}(M)\to C^\infty(M)
\]
is given by
\[
R(X_1,X_2,X_3,X_4)= (R(X_3,X_4)X_2,X_1), \quad X_1,X_2,X_3,X_4\in
{\mathcal X}(M).
\]

If $\dim M=2$, then the Riemannian curvature tensor $R$ has $4$
nontrivial components
\[
R_{1212}=-R_{2112}=R_{1221}=-R_{1221},
\]
Other components equal zero. We have
\[
2R_{1212}=R (g_{11}g_{22}-g^2_{12}),
\]
moreover the scalar curvature $\kappa$ is related with the Gauss
curvature $K$ by
\[
R=2K.
\]

Now assume that $M$ is two-dimensional, and $\gamma$ is a one
dimensional smooth submanifold (a closed curve). Let $s$ be the
natural parameter along $\gamma$. Assume that $\gamma$ is oriented,
and choose the external unit normal vector $N$ at each point of
$\gamma$. In a tubular neighborhood $U$ of $\gamma$, consider the
associated Fermi coordinate system $(X_0,X_1)=(s,t)$, which is
defined as follows. For any $x\in U$, $t$ is the distance from $x$
to $\gamma$ and $s$ is the coordinate of the intersection point of
the minimal geodesic $\xi$, passing through $x$ orthogonally to
$\gamma$. Then at each point $x\in U$ the vector
$N=\frac{\partial}{\partial t}$ coincides with the unit tangent
vector of the minimal geodesic $\xi$, passing through $x$
orthogonally to $\gamma$. Thus, by definition, we have
\[
\|N\|=1, \quad \nabla_NN=0.
\]
Consider the vector field $A$ in $U$ defined in local coordinates
$(X_0,X_1)=(s,t)$ as
\[
A=\frac{\partial}{\partial s}.
\]
For any $x\in \gamma$ the vector $A(x)$ is tangent to $\gamma$ and
\[
\|A(x)\|=1.
\]
Observe that
\[
[A,N]=0.
\]
Then we have
\[
\nabla_N\nabla_NA+R_{AN}N=0.
\]
Indeed, by definition, we have
\[
R_{AN}N=\nabla_A\nabla_NN-\nabla_N\nabla_AN-\nabla_{[A,N]}=-\nabla_N\nabla_AN=-\nabla_N\nabla_NA.
\]
Next, we have
\begin{align*}
N(A,N) & = (\nabla_NA,N)+(A,\nabla_NN)\\
& = (\nabla_NA,N).
\end{align*}
On the other hand, since $\|N\|\equiv 1$, we obtain
\begin{align*}
0=&A(N,N)= 2(\nabla_AN,N).
\end{align*}
Therefore, we have
\[
N(A,N)=0,
\]
that is, $(A,N)$ is constant along the integral curves of $N$, which
are minimal geodesics orthogonal to $\gamma$. Since $(A,N)=0$ on
$\gamma$, it follows that
\[
(A,N)=0.
\]
Thus, the coefficients of the metric have the form
\[
g_{00}(s,t)=(A,A):=a(s,t), \quad g_{01}(s,t)=(A,N) =0, \quad
 g_{11}(s,t) = (N,N) = 1.
\]
Observing that
\[
a(s,0)=1\,,
\]
the metric $g$ has the form
\[
g=a(s,t)ds^2+dt^2.
\]
Let us compute Taylor's expansion of $a$ at $t=0$:
\[
g_{00}(s,t)=1+a_{1}(s)t+{a}_{2}(s)t^2+\mathcal O(t^3).
\]
First, we recall the definition of the  mean curvature of $\gamma$. Since $\|A\|=1$ on
$\gamma$, we have
\[
\kappa=(T_AA,N)=(\nabla_AA,N).
\]
Thus we obtain that
\begin{equation}\label{e:A1}
a_1(s) =N(A,A) = 2(\nabla_NA,A)=
2(\nabla_AN,A)=-2(\nabla_AA,N)=-2\kappa.
\end{equation}
Then,
\begin{align*}
a_2(s)= & \frac12 N^2(A,A)
= (\nabla_N\nabla_NA,A)+ \|\nabla_NA\|^2.
\end{align*}
For the first term, we get
\[
(\nabla_N\nabla_NA,A)=-(R_{AN}N,A)=-R_{1001}=-\frac12R
\]
For the second term, since $(A,N)=0$ on $\gamma$, we obtain
\[
\|\nabla_NA\|^2=(\nabla_NA,A)^2+(\nabla_NA,N)^2=\kappa^2.
\]
Thus, we have
\begin{equation}\label{e:A2}
a_2=-\frac12R+\kappa^2.
\end{equation}

\end{document}